\documentclass{article}

\PassOptionsToPackage{numbers, compress}{natbib}

\usepackage[preprint]{neurips_2021}

\usepackage[utf8]{inputenc} 
\usepackage[T1]{fontenc}    
\usepackage{amsmath}
\usepackage{hyperref}       
\usepackage{url}            
\usepackage{booktabs}       
\usepackage{amsfonts}       
\usepackage{nicefrac}       
\usepackage{microtype}      
\usepackage{graphicx}

\def\cF{{\cal F}}

\def\vx{\mathbf{x}}
\def\vy{\mathbf{y}}
\def\vr{\mathbf{r}}
\def\vs{\mathbf{s}}

\def\RR{\mathbb{R}}

\newcommand{\wbnn}{\texttt{WideBNet}~} 
\newcommand{\nbnn}{\texttt{NarrowBNet}~} 
\newcommand{\unet}{\texttt{U-Net}~}
\newcommand{\switchnet}{\texttt{SwitchNet}~}
\newcommand{\slbnn}{\texttt{SwitchlessBNet}~}

\usepackage{xcolor}

\usepackage{subfigure}

\usepackage{enumitem}

\title{Accurate and Robust Deep Learning Framework for Solving Wave-Based Inverse Problems in the Super-Resolution Regime}

\author{%
  Matthew Li \\ 
  Computational Science and Engineering\\
  Massachusetts Institute of Technology\\
  Cambridge MA 02139 \\
  \texttt{mtcli@mit.edu} \\
   \And
    Laurent Demanet\\ 
  Department of Mathematics\\ 
  Earth Resources Lab\\
  Massachusetts Institute of Technology\\ 
  Cambridge MA 02139 \\
  \texttt{laurent@math.mit.edu} \\
   \And
    Leonardo Zepeda-N\'u\~nez. \\ 
  Department of Mathematics\\
  University of Wisconsin-Madison\\
  Madison WI 53706 \\
  \texttt{zepedanunez@wisc.edu} \\
}

\begin{document}

\maketitle

\begin{abstract}
We propose an end-to-end deep learning framework that comprehensively solves the inverse wave scattering problem across all length scales.
Our framework consists of the newly introduced wide-band butterfly network \cite{MTC_LD_LZN:Wideband_butterfly_net} coupled with a simple training procedure which dynamically injects noise during training. While our trained network provides competitive results in classical imaging regimes, most notably it also succeeds in the super-resolution regime where other comparable methods fail. This encompasses both (i)~reconstruction of scatterers with sub-wavelength geometric features, and (ii)~accurate imaging when two or more scatterers are separated by less than the classical diffraction limit. We demonstrate these properties are retained even in the presence of strong noise and extend to scatterers not previously seen in the training set. In addition, our network is straightforward to train requiring no restarts and has an online runtime that is an order of magnitude faster than optimization-based algorithms. We perform experiments with a variety of wave scattering mediums and we demonstrate that our proposed framework outperforms both classical inversion and competing network architectures that specialize in oscillatory wave scattering data.

\end{abstract}

\section{Introduction}

Wave-based inverse problems\textemdash i.e., inferring the properties of an unknown medium from indirect observations of the medium's response to probing waves\textemdash are ubiquitous in engineering and sciences. This encounters applications in 
a myriad of different fields such as geophysics, astronomy, biomedical imaging, radar, spectrography, signal processing, communications, among many others \cite{Cheney_SAR:2001,Virieux_FWI:2017,Biondi:3D_seismic_imaging,MRI_Schenk:1996,Natterer_CT:2001}. 

All of these applications rely on the physical ability of waves to propagate information across long distances with little distortion. The amount of information that is communicated by waves is essentially proportional to its carrier frequency: higher frequencies transmit more information and thus, when used for imaging, lead to reconstructions with higher resolution. However, in order to reduce computational cost, classical inversion methods \cite{Oldham1906,Gutenberg1914,Backus_Gilbert:1968,Stefanov2019} selectively ignore later scattering events in the time series or process only the low frequency content in the data. These algorithms therefore offer fast reconstruction \cite{ultra_sound_2015,fast_ultra_sound,Javaherian_2020} of relatively low resolution images whose quality deteriorates significantly for highly heterogeneous media.

Modern techniques\footnote{See \cite{Colton_Kress:Inverse_Acoustic_and_Electromagnetic_Scattering_Theory} and \cite{Virieux_Operto:An_overview_of_full-waveform_inversion_in_exploration_Geophysics} for excellent historical reviews.} have since been developed \cite{RTM:1983,Colton_Kirsh:LSL1996,Kirsh:Factorization_methods,Taratola:Inversion_of_seismic_reflection_data_in_the_acoustic_approximation,Pratt:Seismic_waveform_inversion_in_the_frequency_domain;_Part_1_Theory_and_verification_in_a_physical_scale_model,Virieux_FWI:2017} to service applications which depend on accurately imaging medium discontinuities and other fine grained material details. Unfortunately, theses methods rely on high-dimensional, highly non-linear optimization programs which are susceptible to local minima. Thus, in industrial practices, domain experts are required in the loop to monitor and manually adjust the descent path. This exacerbates the already prohibitively expensive computational cost of these algorithms which require state-of-the-art computer clusters to even process the field data.

Evidently, applications requiring real time reconstructions, such as biomedical imaging, radar, continuous monitoring for CO$_2$ sequestration and geothermal energy, etc, remain beyond the reach of such algorithms \cite{Lucka:FWI_3D_breast_imaging}. 

In this context we seek to apply advances in machine learning to develop high-resolution, on-the-fly, reconstruction methods. Broadly, existing research in this area can be categorized as: 
(i) ``physics informed'' networks which minimize a PDE misfit function with a neural network ansatz \cite{PINN_Inverse_Problems}. While these methods offer extremely promising results they cannot yet provide real-time reconstruction. 
Alternatively, (ii) the mapping from the data to medium perturbations can be learned directly using a neural network \cite{MNNH,MNNH2,FanYing:scattering,FanYing:traveltime,Khoo_YingSwitchNet:2019,Yingzhou2018,Barbastathis_dynamical_ML_inversion}.

For methods in the latter category the primary challenge lies in capturing the non-compressible long-range interactions of the inverse data-to-perturbation map efficiently. Indeed the long-range interactions often cause failures when applying typical approaches based on U-nets \cite{U-Net, MNNH, MNNH2}, as well as for generative modelling tools used for imaging such as VAEs or GANs \cite{Isolapix2pix:2016,Ledig2017,mirzacGANs:2014}. Recent approaches \cite{FanYing:scattering,FanYing:traveltime,Khoo_YingSwitchNet:2019,Yingzhou2018} address this issue by embedding the physics of wave propagation into the architecture of the network; specifically, taking cues from Fourier integral operators (FIOs) \cite{Hormander:The_Analysis_of_Linear_Partial_Differential_Operators_FIO} and leveraging the butterfly algorithm \cite{Candes_Demanet_Ying:A_Fast_Butterfly_Algorithm_for_the_Computation_of_Fourier_Integral_Operators,Borm_Butterfly:2017} for fast application of such operators. This builds on the use of FIOs to describe a diverse set of inverse problem modalities such as reflection seismology \cite{Lin_inv_scattering}, thermoacoustic and photoacoustic tomography, radar \cite{FIO_radar,Cheney_SAR:2001,Quinto_2011}, and single photon emission computed tomography \cite{micro_local_Spect}. However, these methods have several drawbacks: (i) their complexity increases super-linearly with respect to the Shannon-Nyquist scaling,  (ii) they depend on stringent assumptions on symmetries and other equi- and/or in- variances, or (iii) they exhibit unstable training dynamics due to the highly-oscillatory nature of the data \cite{F_principle} and require intrusive initializations to obtain accurate approximations \cite{Yingzhou2018}. 

The recently introduced wide-band butterfly network (\texttt{WideBNet}) addresses these pitfalls in a unified manner \cite{MTC_LD_LZN:Wideband_butterfly_net}. Their main contribution was in extending the butterfly networks \cite{Butterfly-Net2,Yingzhou2018} to assimilate wide-band data in a multi-scale fashion similar to the Cooley-Tukey algorithm \cite{Cooley_Tukey:1965}. This fully realizes the advantages of the butterfly algorithm in enabling a network design with training parameters and complexity which scales linearly, up-to-poly-log factors, with the dimension of the data. Additionally, the use of multi-frequency data stabilizes the inverse problem, and therefore the training dynamics of the network, and also empirically enables imaging in the super-resolution regime. 

Despite the advantages of the \wbnn documented in \cite{MTC_LD_LZN:Wideband_butterfly_net} the stability of this network with noise remains a key gap that remains unaddressed. This is important as the mathematical super-resolution problem is known to be poorly conditioned~\cite{Donoho:super_resolution1992}. Demonstrating robustness against noisy data is therefore key for applications, where the data acquisition is often noisy due to constraints in the exact position of the receivers or due to intrinsic statistical uncertainty. 

\noindent \textbf{Our contribution:} We introduce the first deep learning framework capable of successfully solving the inverse scattering problem in the super-resolution regime while being robust to noise. The resulting network is able to handle multiscale features, with both sub-wavelength and wavelength level features, as well as resolve colliding scatterers separated by less than the diffraction limit. Remarkably, both of these properties hold with resilience to noise -- we demonstrate robustness even with 100\% noise-to-signal corruptions in the data. In addition, our network outperforms both classical and newly introduced architectures such as U-net (\texttt{U-Net}) \cite{U-Net}, SwitchNet (\texttt{SwitchNet}) \cite{Khoo_YingSwitchNet:2019} and butterfly-networks (\texttt{NarrowBNet}) \cite{Butterfly-Net2,Yingzhou2018} for a comprehensive range of scatterers with varying wave-scattering physics. To our knowledge this is the first network to achieve these feats, thus enabling future applications in biomedical imaging, radar, and geophysical monitoring, where real-time imaging across all length scales is mission critical.

\section{Physical Model}
\label{sec:phys_model}
\begin{figure}
    \centering
    \includegraphics[width=0.32\textwidth]{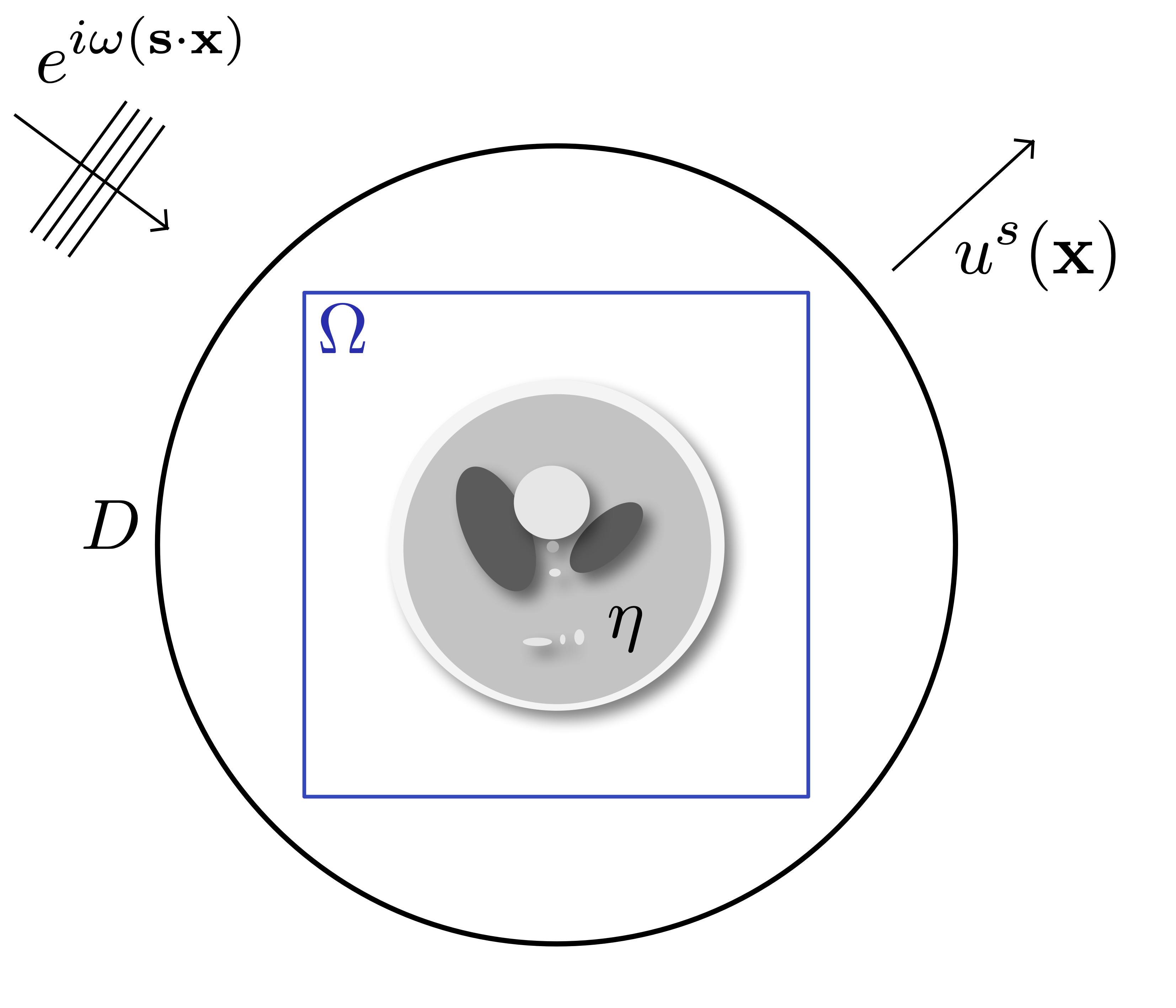}
    \includegraphics[width=0.32\textwidth]{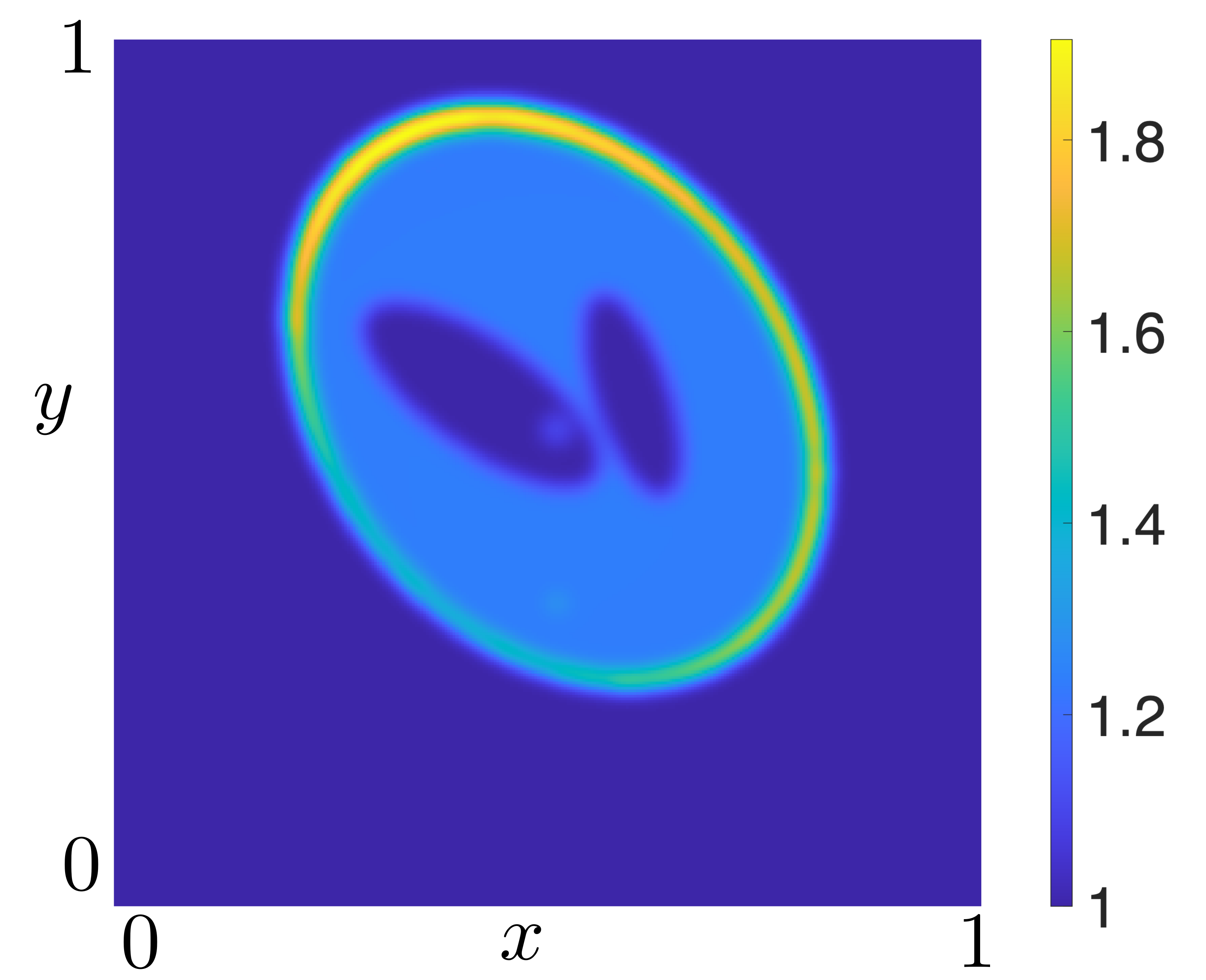}
    \includegraphics[width=0.32\textwidth]{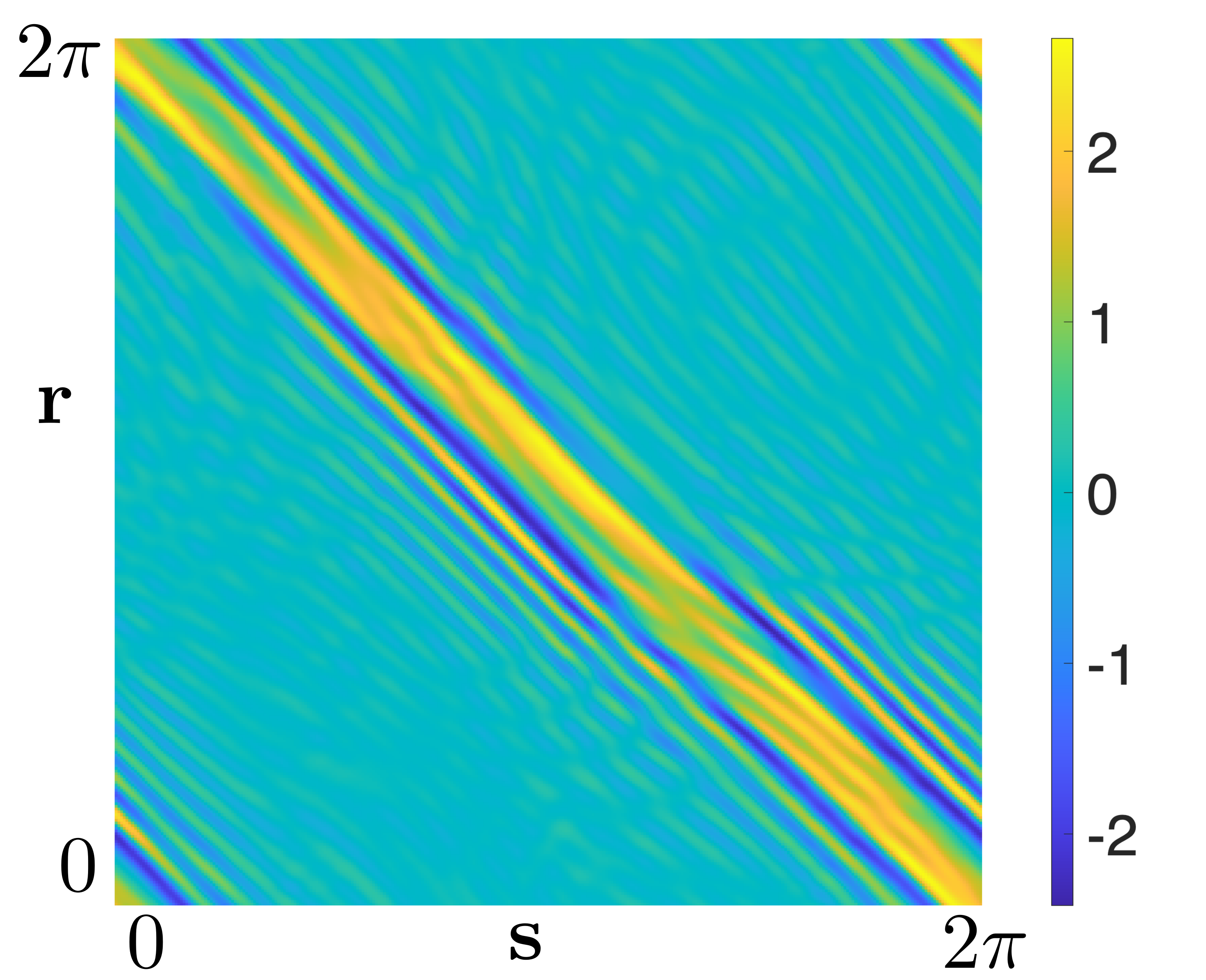}
    \caption{{\it (left)} setup for the inverse scattering problem. In the figure we probe the medium with a planewave with direction $\vs$, and we sample the scattered field on the disk $D$, {\it (center)} slowness squared of a Shepp-logan model and {\it (right)} real part of the discrete Dirichlet-to-Neumann map.
    }
    \label{fig:scattering}\textbf{}
\end{figure}

We consider the time-harmonic wave equation with constant-density acoustic physics, also called the Helmholtz equation, with frequency $\omega$ and squared slowness $m$, given by  
\begin{equation} \label{eq:helmholtz}
    (\Delta + \omega^2 m(\vx)) u (\vx) = 0
\end{equation}
with radiating boundary conditions. We further suppose the slowness admits a scale separation into
\begin{equation}
    m(\vx) = m_0(\vx) + \eta(\vx) = 1 + \eta(\vx),
\end{equation}
where $m_0$ corresponds to the smooth background slowness, assumed known and for simplicity equal to one\footnote{This assumption is only made to make the presentation more transparent.}, and $\eta$ the rough perturbation that we wish to recover. Under these assumptions, solutions to \eqref{eq:helmholtz} can be expressed in the form
\begin{equation}
    u(\vx) = e^{i\omega \left ( \mathbf{s} \cdot \vx \right )} + u^{s}(\vx), 
\end{equation}
where $e^{i\omega \left ( \mathbf{s} \cdot \vx \right ) }$ is the incoming plane wave, with propagating direction $\mathbf{s}$, that we use to ``probe'' the perturbation, and $u^{s}(\vx)$ is the scattered field produced by the interaction of the perturbation with the impinging wave. The scattered field satisfies \cite{Colton_Kress:Inverse_Acoustic_and_Electromagnetic_Scattering_Theory}
\begin{equation} \label{eq:scattering}
    \left \{   \begin{array}{ll}
                \left(\Delta+\omega^{2}( 1 + \eta(\vx)\right) u^{s}(\vx) = -\omega^2 \eta(\vx) e^{i \omega ( \mathbf{s} \cdot \vx)} & \textrm{ for }\,\, \vx \in \mathbb{R}^2,  \\ \displaystyle
                \lim _{|\vx| \rightarrow \infty}|\vx|^{1/2}\left(\frac{\partial}{\partial|\vx|}-\mathrm{i} \omega\right) u^{\mathrm{s}}(\vx)=0,
                \end{array}
    \right .
\end{equation} 
following the setup depicted in Figure~\ref{fig:scattering}.
For simplicity, we select the detector manifold~$D$ to be a circle of radius~$R$ that surround the domain of interest $\Omega$ as shown in Figure~\ref{fig:scattering} {\it (left)}. For each incoming direction $\mathbf{s} \in \mathbb{S}^1$ the data is given by sampling the scattered field with receiver elements that are located on $D$ and indexed by $\mathbf{r} \in \mathbb{S}^1$. This yields the far-field pattern given by $\Lambda_{\vs, \vr}^{\omega} = u^{sc}(R\vr)$, where $\vs$ denotes the incoming probing direction as defined in \eqref{eq:scattering}. We call $\mathcal{F}^{\omega}[\eta]:= \Lambda_{\vs, \vr}^{\omega}$ the {\it forward map} relating the perturbation $\eta$ to its corresponding far-field pattern. 
In Figure~\ref{fig:scattering} {\it (right)} we can observe the oscillatory behavior of a typical example of the far-field pattern for the perturbation in Figure~\ref{fig:scattering} {\it (center)}. 

Following classical WKBJ asymptotic theory \cite{olver_2014:asymptotics} and Fourier analysis \cite{Hormander:The_Analysis_of_Linear_Partial_Differential_Operators_FIO} one can write 
\begin{equation} \label{eq:Fourier_operator}
 \Lambda_{\vs, \vr}^{\omega} = (\mathcal{F}^{\omega} \eta)(\vx)=\int_{\text{supp}(\eta)} a(\vx, \vy) e^{i \omega \phi(\vx, \vy)}  \eta(\vy) d \vy.
\end{equation}
where $\vx = (\vs, \vr)$, $a$ and $\phi$ are functions. In this case $a$ and $\phi$ are solutions to a transport and an eikonal equation respectively, that depend on $\eta$. Fortunately, in the Born approximation this operator can be reduced to a Fourier transform (see \eqref{eq:far_field_pattern} in Appendix \ref{app:FBP}). 
 
We define the \emph{inverse map} as the map estimating $\eta$ from $\Lambda^{\omega}_{\vs,\vr}$, namely $\left ( \cF^{\omega} \right)^{-1}$. The inverse problem can be formulated as a minimization problem, i.e., 
\begin{equation} \label{eq:fwi}
   \left ( \cF^{\omega} \right)^{-1} \left [\Lambda_{\mathbf{s},\mathbf{r}}^{\omega} \right] =  \eta^* := \text{argmin}_{\mu} \|\mathcal{F}^{\omega}[\mu] - \Lambda_{\mathbf{s},\mathbf{r}}^{\omega} \|^2
\end{equation}
for a suitable set of perturbations, and for a suitable norm. This problem is often solved using PDE-constrained optimization methods, or in the linearized case, using filtered back-projection (see Appendix \ref{app:FBP}). However, this has been proven to result in unstable algorithms which can converge to non-physical local minima \cite{Hahner_Hohage2001_inverse_problem_estimates}. Thus, one key insight for designing robust algorithms, which is also incorporated in wide-band butterfly network, is to use wide-band data viz.,
\begin{equation} \label{eq:fwi_weighted}
    (\mathcal{F})^{-1}[\Lambda_{\mathbf{s},\mathbf{r}}^{\omega_1}, \Lambda_{\mathbf{s},\mathbf{r}}^{\omega_2}, ..., \Lambda_{\mathbf{s},\mathbf{r}}^{\omega_{n_\omega}}] = \eta^* := \text{argmin}_{\mu} \sum_{i = 1}^{n_{\omega}} w_i\|\mathcal{F}^{\omega_i}[\mu] - \Lambda_{\mathbf{s},\mathbf{r}}^{\omega_i} \|^2.
\end{equation}
How to choose the frequencies $\{\omega_i\}_{i = 1}^{n_{\omega}}$ and the corresponding weights $\{w_i\}_{i = 1}^{n_{\omega}}$ remains an area of research \cite{Borges_Gillman_Greengard:2017}. In fact, this problem still present spurious local minima, thus regularization techniques such as recursive linearization are used \cite{Chen:Inverse_scattering_via_Heisenberg's_uncertainty_principle,Pratt:Seismic_waveform_inversion_in_the_frequency_domain;_Part_1_Theory_and_verification_in_a_physical_scale_model}. In a nutshell, recursive linearization (or frequency sweeps) can be understood at the algorithmic level as weights that change throughout the optimization loop: At the beginning only the low-frequency weights are not zero, and  during the optimization the mass of the weights is slowly transferred to the weights involving higher frequencies. This is another key insight incorporated in to the design of the wide-band butterfly network: Instead of treating each frequency independently, processing low-frequency data provides guidance to process the higher-frequency data.

For the discrete problem we assume that $\Omega$ is discretized in $n_x \times n_y$ equispaced points in space, and that the far field data is probed following $n_{sr}$ equispaced direction in $\mathbb{S}$, and that the scattered field is probed at $n_{rs}$ equispaced points in $D$. For simplicity we assume that $n_x = n_y = n_{sr} = n_{rs}$.

\section{Wideband Butterfly Network Architecture}

\begin{figure}
    \centering
    \includegraphics[trim = 25mm 9mm 25mm 9mm, clip, width=0.98\textwidth]{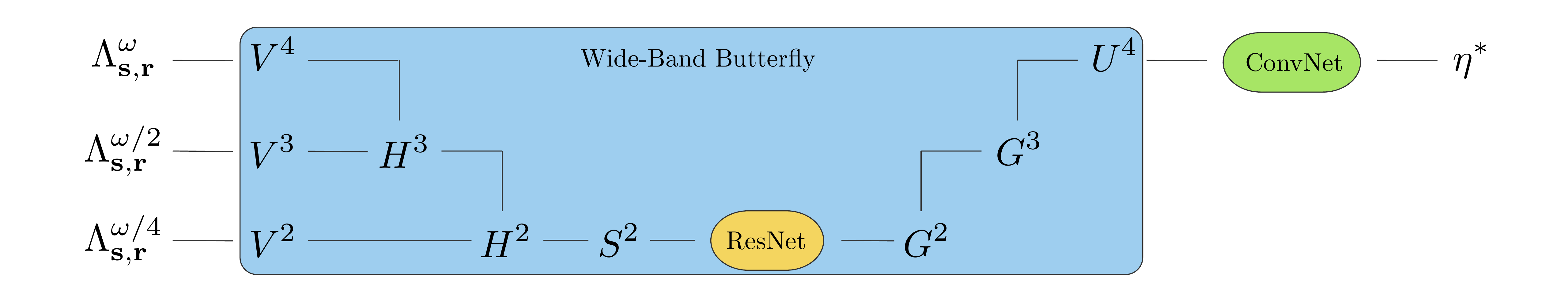}
    \caption{Sketch of the architecture of the Wide-band buttefly network for $L =4$ used to approximate $\eta^* = (\mathcal{F})^{-1}[\Lambda_{\mathbf{s},\mathbf{r}}^{\omega}, \Lambda_{\mathbf{s},\mathbf{r}}^{\omega/2}, \Lambda_{\mathbf{s},\mathbf{r}}^{\omega/4}]$ .}
    \label{fig:WBNet}
\end{figure}

The network is composed of a convolutional block and a wide-band butterfly block with a non-linear resnet in the middle following Figure~\ref{fig:WBNet}. This construction is inspired by the filtered back-projection, which is ubiquitous in inverse problems (see Appendix \ref{app:FBP}), and the butterfly factorization \cite{Butterfly-factorization:Liu_Xing2020,Li_Yang_Martin_Ho_Ying:Butterfly_Factorization}. In a nutshell, \wbnn \emph{lifts} the linearized version of the inverse problem (see  \eqref{eq:FBP} in Appendix \ref{app:FBP}) with the butterfly algorithm to a non-linear operator that is fully trainable. The inputs to the network are the data $\Lambda^{\omega}_{\vs, \vr}$ at frequency $\omega$, each of which is an $n_{sr} \times n_{rs}$ matrix reshaped into a tensor of dimensions  $(2^L, 2^L, s^2)$. The first two dimensions of the reshaped tensor account for the geometric information in $2^L \times 2^L$ cells, and $s^2$ corresponds the size of the leaf nodes (see Figure~\ref{fig:input_ordering}).

Specific to the \wbnn are the local embedding blocks $V^{\ell}$, the down-sampling layers $H^{\ell}$, the switch layer $S^{L/2}$, the upsampling layers $G^{\ell}$, and the interpolation layer $U^{L}$. For simplicity we describe the action of these components at the highest frequency, i.e. at the highest length scale $L$. For the full description we refer to the Appendix~\ref{app:architectures} or the original manuscript \cite{MTC_LD_LZN:Wideband_butterfly_net}.  The implementation of the convolutional and resnet layers are standard and details on the dimensions are similarly left to the Appendices \ref{app:train_details} and \ref{app:architectures}.

\begin{itemize}[leftmargin=*]
    \item{\textbf{Local embedding layer} $\mathbf{V^{\ell}}$}  takes input dimension  $(2^{\ell}, 2^{\ell}, s^2)$ and processes each cell to produce a local representation as a vector of dimension $r$, thus outputting a tensor of dimensions  $(2^{\ell}, 2^{\ell}, r)$.  The first two components index the geometrical position of each cell and the last corresponds to the number of channels in the local embedding (see Figure~\ref{fig:V_ell} for a sketch of this operation).
    
    \item{\textbf{Downsampling layer} $\mathbf{H^{\ell}}$}  merges information from contiguous cells, They downsample the number of cells by two in each axis and increase the information contained in the new aggregate cell by a factor four. If the input is $(2^\ell,2^\ell, r)$ then the output is $(2^{\ell-1},2^{\ell-1}, 4r)$ (see Figure~\ref{fig:H_ell} for a sketch of this operation).
    
    \item{\textbf{Switch layer} $\mathbf{S^{L/2}}$} was introduced in \cite{Khoo_YingSwitchNet:2019} to perform an all-to-all communication between the cells. The input and output dimension $(2^{L/2},2^{L/2}, 2^{L}r)$ since this layer broadcasts information in the third dimension across the first two in packages of size $r^2$ (see Figure~\ref{fig:switch_ordering} for a sketch of this operation).
    
    \item{\textbf{Upsampling layer} $\mathbf{G^{\ell}}$} is formally the adjoint of the downsampling layers, i.e., they split information in each cells into four contiguous cells, each with a fourth of the original aggregate cell. If the input is $(2^\ell,2^\ell, 4r)$ then the output is $(2^{\ell+1},2^{\ell+1}, r)$.
    
    \item{\textbf{Local sampling layer} $\mathbf{U^{L}}$} takes each cell of the input of dimension  $(2^L, 2^L, r)$ which contains a local representation and they sample it in the corresponding leaf of size $s^2$, thus the output is  $(2^L, 2^L, s^2)$. 

\end{itemize}

Observe from Figure~\ref{fig:WBNet} that the information is down-sampled after each $H^{\ell}$ layer to inject information at the correct length scales. This complies with Shannon-Nyquist sampling theory which describes the relation of the information content and sampling frequency.

\section{Numerical Results} \label{sec:numerical_results}

We benchmark the \wbnn architecture for inversion of scatterers with various geometric configurations and length scales. Our imaging domain spans $[-0.5, 0.5]^{\otimes 2}$ and is discretized over an equispaced mesh with $n = 80$ points per dimension. This corresponds to a quad-tree partitioning with $L=4$ levels and $s=5$ leaf nodes. We consider $n_{\omega} = 3$ number of frequencies at $2.5$, $5.0$, and $10.0$ Hz in which the data was assimiliated at levels $\ell=2$,~$3$, and $4$ respectively. We consider a homogeneous background squared slowness of $m_0 = 1$ (or equivalently a background wave-speed of $c_0 = 1$) which corresponds to $8$~points per wavelength (PPW) at the highest frequency and limits the classical Rayleigh resolution to $\lambda_\textrm{min}/2 = 4$~pixels. See Appendices \ref{app:train_details} and \ref{app:architectures} for details the choice of hyperparameters and implementation details, respectively.

The synthetic scattering data was generated with second-order finite differences for training and fourth-order finite differences for testing. The use of different stencils validates against overfitting of numerical dispersion artifacts. The radiating boundary conditions in \eqref{eq:scattering} were implemented using perfectly matched layers (PML) spanning one wavelength with quadratic profile of intensity~$80$~\cite{Berenger1994}. The scattered waves were sampled at $80$~equi-angular intervals on a circle of radius $r=0.5$ where the incident waves arrived at angles aligned with the receiver geometry. The training and testing split was $21000$ points and $4000$ points, respectively. Note each datapoint is a tensor of dimension~$n_{sr} \times n_{rs} \times n_{\omega}$, where $n_{sr} = n_{rs} = 80$ and $n_{\omega} = 3$.

The data at each frequency was normalized by subtracting and scaling by the mean and variance of the aggregate pixel intensities in the training dataset. Whereas \cite{MTC_LD_LZN:Wideband_butterfly_net} examines \wbnn with noiseless data, here we consider a multiplicative noise model similar to \cite{Khoo_YingSwitchNet:2019} with 
\begin{equation} \label{eq:noise:model}
\Lambda_{\vs,\vr} = \Lambda_{\vs,\vr}^{\text{clean}}(1 + \epsilon(\vs,\vr)),\quad \epsilon(\vs,\vr)\stackrel{iid}{\sim} \mathcal{N}(0, \sigma^2).
\end{equation}
Since each frequency is standardized independently this corresponds to multiplicative coloured noise in the time series data. We set $\sigma=1$ to introduce $100\%$ noise to signal ratio in our data. The noise is applied to each data pixel dynamically each epoch. 

The objective function we used was the pixelwise sample loss 
\[
\sum_{x = \text{pixel in image}} \lVert (K_{\text{high}} \ast \eta)(x) - \wbnn[\Lambda_{s,r}^{\omega_1}, \ldots, \Lambda_{s,r}^{\omega_{n_{\omega}}}](x)  \rVert_2^2
\]
where $\eta$ is the total scattering medium and $\{\Lambda_{s,r}^{\omega_{i}}\}_{i= 1}^{n_{\omega}}$ the noisy multifrequency data. The scattering medium data was smoothed using a high-pass filter $K_\text{high}$ described by a Gaussian kernel with characteristic width of $0.75$ grid points. This was observed to promote faster training in~\cite{MTC_LD_LZN:Wideband_butterfly_net} for noiseless data. Since the support of this filter is significantly less than the Nyquist limit this processing does not pollute or eliminate sub-wavelength features in the medium. 

For completeness we also investigate the performance of comparable architectures such as: the popular image processing architecture \unet  
where we feed in multifrequency data~\cite{U-Net};  and 
\texttt{SwitchNet}, which similarly uses the butterfly architecture but requires superlinear $\mathcal{O}(N^{3/2})$ complexity for training stability when using single frequency data. We also consider a narrow-band implementation of our \wbnn architecture referred to as 
\nbnn which uses only $10$~Hz data assimilated at level $L=4$ and otherwise coincides with the architecture proposed in \cite{Yingzhou2018}. Furthermore we benchmark 
\texttt{SwitchlessBNet}, which has identical degrees to freedom as \wbnn except for the omission of the $S^{L/2}$ layer. We refer the interested reader to Appendices \ref{app:FBP} and \ref{app:architectures} for a comparison against more classical methods such as filtered back-projection (see Figure \ref{fig:FIO_wbnet} and PDE-constrained optimization methods (see Figure \ref{fig:fwi_bdnet_comparison}), in which \wbnn outperforms both methods, particularly in the noisy case. 

\subsection{Performance across classes}
 
\begin{figure}
    \centering
    \includegraphics[width=0.99\textwidth]{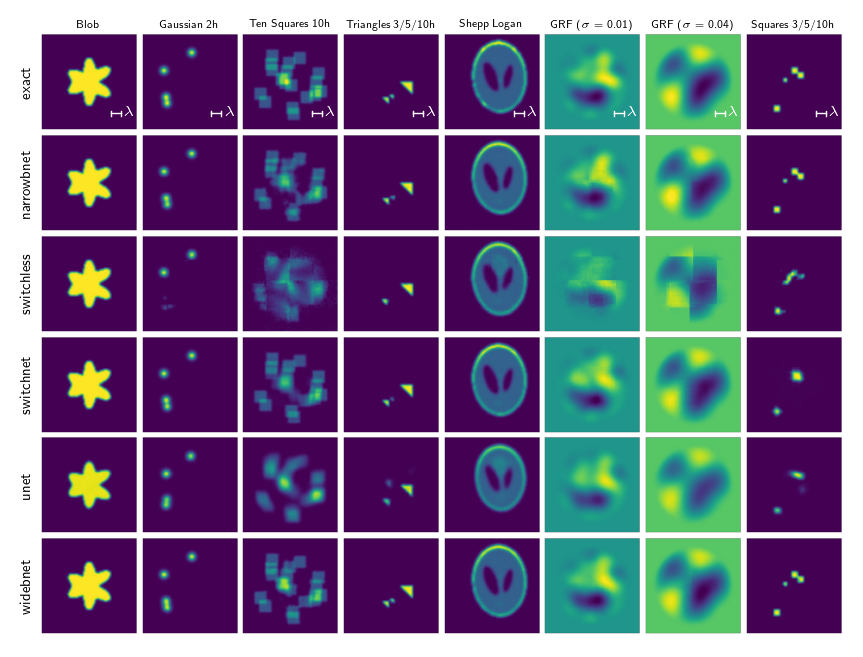}
    \caption{Comparison of performance between different architectures (rows) on various datasets (columns). The scattering data was contaminated with $\sigma=1$ multiplicative noise following \eqref{eq:noise:model}.}
    \label{fig:testing}
\end{figure}

Table~\ref{tab:benchmark} and Figure~\ref{fig:testing} depicts the performance of various architectures for different scattering media. These media comprehensively cover a wide range of wave scattering behaviour: \textbf{Randomized Shepp-Logan} is a standard biomedical imaging dataset which emulates contrasting material properties inside the human head \cite{random_shepp_logan}; \textbf{Gaussian Random Field} propagate waves with strong multipathing and weak coda. The parameter $\sigma = \{0.01, 0.04\}$ refers to the characteristic length scale of the fluctuations; \textbf{Blob}  contains a single non-convex scatterer with wavelength level features, \textbf{Ten-Squares 10h} corresponds to ten randomly located squares each with 10 pixel sidelengths. Waves traveling in this medium exhibit strong multiple scattering; \textbf{Gaussian 2h} correspond to point scatterers at the highest resolution limits of the medium; \textbf{Squares 3/5/10h} are randomly located squares of $3$, $5$, and $10$ pixel sidelengths. This covers multiple lengthscales ranging from subwavelength, wavelength, and superwavelength scatterers; \textbf{Triangles 3/5/10h} are similar to squares except with right triangles. For the geometric scatterers the amplitude of each scatterer is fixed at $\Delta m = 0.2$, however superposition of scatterers can lead to higher local wavespeeds. The geometric mediums contain uniformly sampled range of 2, 3, or 4 scatterers which are uniformly located in a circle of radius $r=0.35$. All scattering mediums lie in a homogeneous background of $m_0 =1$. 

\begin{table}[]
\centering
\caption{Mean Squared Error per pixel on testing data with $\sigma=1$ multiplicative noise.}
\label{tab:benchmark}
\begin{tabular}{@{}lccccc@{}}
\toprule
 & \multicolumn{1}{c}{WideBNet} & \multicolumn{1}{c}{UNet} & \multicolumn{1}{c}{SwitchlessBNet} & \multicolumn{1}{c}{NarrowBNet} & \multicolumn{1}{c}{SwitchNet} \\ \midrule
Blob & \textbf{5.5996E-04} & 4.3213E-03 & 2.6254E-03 & 5.9252E-04 & 6.9086E-04 \\
Square 3/5/10h & \textbf{7.9643E-03} & 1.3671E-01 & 1.5380E-01 & 7.9946E-03 & 1.9018E-01 \\
Triangle 3/5/10h & 8.1558E-03 & 1.4777E-01 & 1.3805E-01 & 8.7101E-03 & \textbf{6.3816E-03} \\
Gaussian 2h & 2.5700E-03 & 6.2389E-02 & 1.1298E-01 & 4.7489E-03 & \textbf{2.4208E-03} \\
Ten Squares 10h & \textbf{3.1379E-02} & 1.0738E-01 & 4.0320E-01 & 9.4893E-02 & 8.9998E-02 \\
Shepp Logan & \textbf{7.2137E-03} & 3.8251E-02 & 3.1400E-02 & 1.1412E-02 & 8.6461E-03 \\
GRF ($\sigma$=0.04) & \textbf{3.3387E-03} & 1.2358E-02 & 2.1289E-01 & 9.0839E-03 & 9.1764E-03 \\
GRF ($\sigma$ = 0.01) & \textbf{1.3585E-02} & 4.4441E-02 & 4.7032E-01 & 5.3533E-02 & 5.7688E-02 \\ \bottomrule
\end{tabular}
\end{table}

Quantitatively in Tab.~\ref{tab:benchmark} we observe that \wbnn generally provides the lowest average testing errors across all datasets. Notably, 
\texttt{SwitchlessBNet}, which has an identical architecture to \wbnn (as seen in Figure~\ref{fig:WBNet}) except eliminating the switch layer, performs two orders of magnitude worse on average. This is consistent with intuition from FIO theory where the switch permutation is necessary to capture the long-range (i.e. local-to-global) interactions inherent to wave scattering data. Although \switchnet also uses this switch permutation layer it does not outperform \wbnn despite requiring an additional million trainable weights in its architecture (see Appendix \ref{app:hyper_parameters} for the total number of trainable parameters of each network). The \unet architecture performs similarly poorly indicating that the multiscale behaviour in the dataset is difficult to capture with purely convolutional filters. Using solely narrow-band (single frequency) data provides comparable behaviour on simple geometric datasets, but suffers when capturing more complex media such as the Shepp-Logan or gaussian random fields. We observe that the Blob dataset is the easiest to invert for all architectures.

Qualitatively in Figure~\ref{fig:testing} we observe that \slbnn is unable to fully resolve all of the gaussian scatterers. Moreover, although it is able to resolve the large length scale scatterers in both the square and triangle datasets, the subwavelength scatterers are notably omitted. These observations indicate that \slbnn struggles in the super-resolution regime. Note that the \unet incorrectly orients the colliding point scatterers for the gaussian dataset. A similar phenomeon for the \unet is observed for the Ten-Squares dataset which features multiple scatterers with subwavelength separation. Both \nbnn and \switchnet also struggle with densely overlapping scatterers, whereas \wbnn does not. This reinforces the claim that wideband data stabilizes inversion. Also observe that \nbnn introduces mild artifacts in the gaussian random field with $\sigma=0.04$, while \switchnet images the wrong square.  In all instances the Shepp-Logan phantom is well resolved with the exception of \unet which incorrectly scales the phantom. Lastly we note that both the \nbnn and \switchnet architecture are surprisingly able to resolve certain subwavelength features despite only provided single frequency data. We speculate that they utilize geometric triangulation to achieve these images.

\subsubsection{Out of Distribution Generalization}
\begin{figure}
    \centering
    \includegraphics[width=0.95\textwidth]{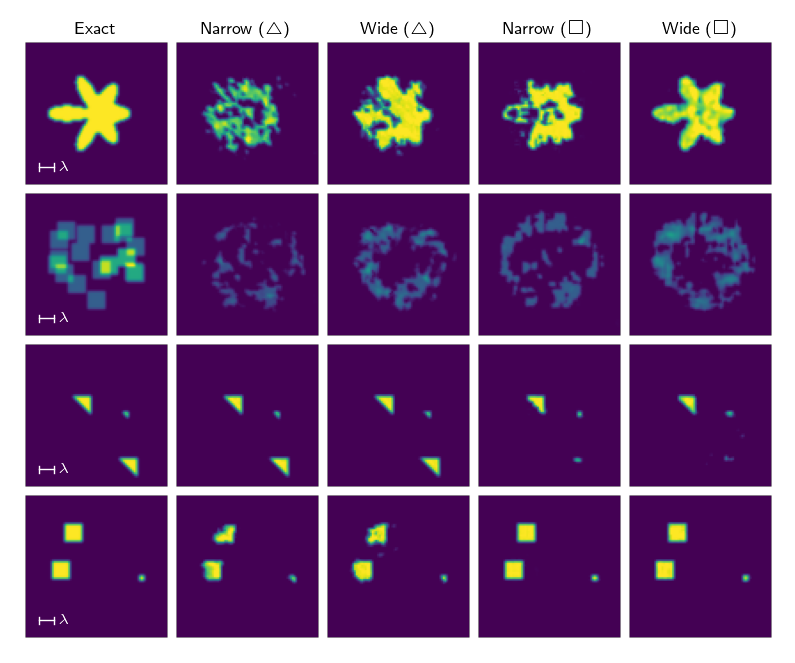}
    \caption{Example of out-of-distribution generalization properties of \nbnn and \texttt{WideBNet}. The shape in parentheses 
    indicate the training set -- $\triangle$ (resp. $\square$) corresponds to Triangle 3/5/10 (resp. Square 3/5/10). 
    Scattering data was polluted by $\sigma=1$ even for the out-of-distribution test cases.}
    \label{fig:generalization}
\end{figure}

We report on the unexpected ability of our trained networks to generalize \textbf{across} datasets. However we caution that this behaviour appears to be strongly dependent on the dataset -- in particular, we only observe successful instances of generalization when training with the Square 3/5/10h and Triangle 3/5/10h datasets. Figure~\ref{fig:generalization} documents this for the narrow- and wide-band architectures where the shape in parentheses in the column titles indicates the training set. A rough numerical comparison of extrapolative properties between datasets can be found in the heatmap shown in Figure~\ref{fig:generalization_heat_map} in the appendix.

Training \wbnn with the square 3/5/10 dataset appears to provide the best out of distribution performance, both qualitatively and quantitatively. For instance, observe that the wide-band network successfully captures the boundaries of the out-of-distribution blob in Figure~\ref{fig:generalization}. We emphasize that training points in the square dataset contain at most four scatterers -- therefore this ability the infill the blob with multiple, i.e. greater than four, scatterers suggests the network has learned some generalized physics of wave scattering. Unexpectedly, however, this behaviour is not repeated for the Ten Squares dataset despite the fact that squares with 10 pixel sidelengths are already present in the training set.

Furthermore, note that both the squared trained \nbnn and \wbnn are able to localize and resolve triangular scatterers. In both instances however the  large triangular scatterer on the bottom is mistakenly characterized as a small triangular scatter; it is unclear why both networks share this deficiency. Surprisingly, triangular trained networks do not generalize as well to the square dataset even though each square can be expressed as a union of two (appropriately rotated) triangles. We leave to future work a systematic investigation of how training datasets can be designed to promote better out of distribution behaviour on trained \wbnn networks. However we note that these results strongly indicate that the networks learn a far more complicated signal processing algorithm than simply joint blind deconvolution with super resolution.

\subsection{Super Resolution beneath the Diffraction Limit}
\label{ssec:superres}
\begin{figure}[h!]
    \centering
    \includegraphics[width=0.98\textwidth]{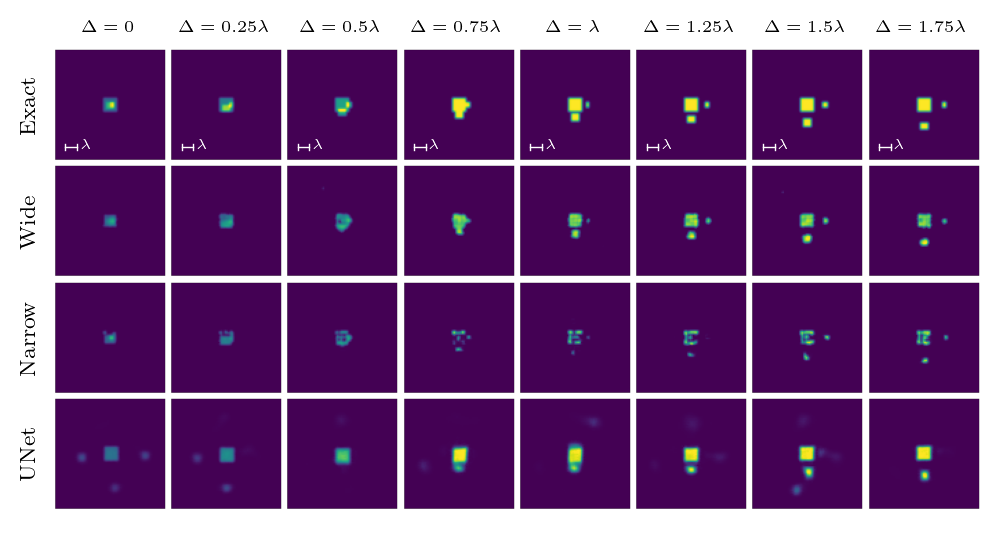}
    \caption{Comparison of imaging capabilities of different networks for square scatterers with side lengths 3, 5, and 10 pixels, which are separated by~$\Delta$. We vary $\Delta$ between sub-wavelength and wavelength separations (one wavelength~$\lambda$ = 8 pixels). }
    \label{fig:collide-square-10-3-5h}
\end{figure}

We report on the claimed super-resolution capabilities of our trained networks: for a comparison with optimization-based methods please see Figure \ref{fig:fwi_bdnet_comparison}. Figure~\ref{fig:collide-square-10-3-5h} depicts the testing performance on a scattering medium with three squares of 3, 5, and 10 pixel side length that are separated by a distance of $\Delta$\footnote{Due to the way the scatterers were defined in code the $\Delta$ is not exactly the separation between the \emph{boundaries} of the scatterers.}. All datasets are polluted with $\sigma=1$ multiplicative noise. When $\Delta \gg \lambda$ we observe that the both \nbnn and \wbnn successfully image all scatterers, although \nbnn struggles with the boundaries of the largest scatterer. However, when $\Delta \leq \lambda$ we note that only \wbnn succeeds in imaging, as claimed. \unet performs poorly over the entire range of $\Delta$ separation as it resolves the largest scatterer but ignores either one or both of the remaining small scatterers.

\subsection{Stability with Noise} \label{sec:noise_stability}
\begin{figure}[h!]
    \centering
    \includegraphics[trim = 2mm 0mm 2mm 0mm, clip, width=0.99\textwidth]{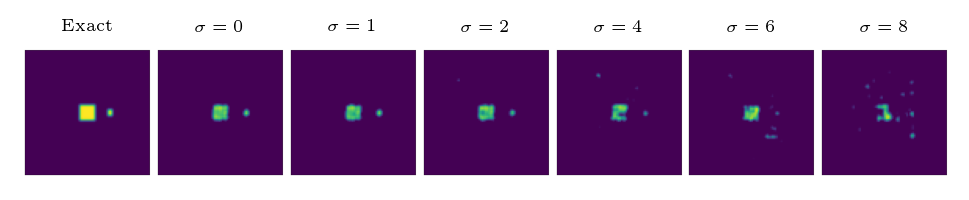}
    \caption{Degradation of image reconstruction with increasing noise-to-signal ratio $\sigma$. The network was trained with data contaminated with $\sigma=1$ noise.}
    \label{fig:noise-sigma}
\end{figure}

Lastly we consider the robustness of the learned networks to noise. Figure~\ref{fig:noise-sigma} depicts the performance of a network trained with $\sigma=1$ noise with testing data that is polluted by different noise levels. The network demonstrates strong denoising abilities as it produces a similar images even for $\sigma = 2$, i.e. a noise level two times higher than what it saw in training. Conversely, the imaging performance also does not suffer when the input data is noiseless with $\sigma=0$. We note that the image quality degrades gracefully with increasing  $\sigma$; while remnants of the large scatterer can still be identified, the small scatterer is not longer reliably imaged at higher noise levels.

\section{Conclusions} \label{sec:conclusions}
We compared \wbnn against a range of competing network architectures on a suite of benchmark scattering configurations. Both quantitative, and qualitatively, \wbnn provides better performance. In particular, \wbnn excels in the super-resolution regime whereas other trained networks fail to consistently produce acceptable images. This super-resolution behaviour was found to extend to higher levels of noise than in the training set, though this breaks down for sufficiently high noise levels. Additional mechanisms for robustness to noise remains an open problem. The \wbnn networks were also found to have surprising generalizability across datasets. However, this behaviour was limited only to certain training datasets. The characterization of ``good'' training datasets, and/or strategies to promote transfer learning across datasets, is left to future work.

\section*{Acknowledgements}
The authors thank Total S.A. for support. L.D. is also supported by AFOSR grant FA9550-17-1-0316. L.Z.-N. is also supported in part by the Wisconsin Alumni Research Foundation, the National Science Foundation under the grant DMS-2012292, and NSF TRIPODS award 1740707.

\bibliographystyle{plain} 
\bibliography{references.bib}

\appendix

\section{Filtered Back-Projection} \label{app:FBP}

We can cast the inverse problem for recovering the perturbation as
\begin{equation} \label{eq:fwi_app}
    \eta^* = \textrm{argmin}_{\mu} \|\mathcal{F}^{\omega}[\mu] - \Lambda_{\mathbf{s},\mathbf{r}} \|
\end{equation}
where $\Lambda_{\mathbf{s},\mathbf{r}}$ is the measured data.  We linearize $\mathcal{F}^{\omega}$ to shed light on the essential difficulties of this problem. Using the classical Born approximation in \eqref{eq:scattering} we obtain that 
\begin{equation}
    u^{sc}(\mathbf{x}) = \omega^2 \int_{\RR^2} \Phi^{\omega}(\mathbf{x}, \mathbf{y})  \eta(\mathbf{y}) e^{i \omega ( \mathbf{s} \cdot \mathbf{y})} d\vy,
\end{equation}
where  $\Phi^{\omega}$ is the Green's function of the two-dimensional Helmholtz equation in homogeneous media, i.e.,  $\Phi^{\omega}$ satisfies
\begin{equation} \label{eq:GreenFunction}
    \left \{   \begin{array}{ll} \displaystyle
                \left(\Delta+\omega^{2}\right) \Phi^{\omega}(\vx, \vy) = -\delta(\vx, \vy) & \textrm{ for }\, \vx \in \mathbb{R}^2,  \\ \displaystyle
                \lim _{|\vx| \rightarrow \infty}|\vx|^{1/2}\left(\frac{\partial}{\partial|\vx|}-\mathrm{i} \omega\right) \Phi^{\omega}(\vx, \vy) =0.
                \end{array}
    \right .
\end{equation} 
Furthermore, we can use the classical far-field asymptotics of the Green's function to express
\begin{equation}
    u^{sc}(R \vr) = -\omega^2 \frac{e^{i \omega R}}{\sqrt{R}} \int_{\RR^2}  \eta(y) e^{i \omega ( \mathbf{s} - \vr) \cdot \mathbf{y}} d\vy + \mathcal{O}(R^{-3/2}).
\end{equation}
Thus, up to a re-scaling and a phase change, the far-field pattern defined in \eqref{eq:far_field_pattern} can be approximately written as a Fourier transform of the perturbation, i.e., 
\begin{equation} \label{eq:far_field_pattern}
    \Lambda_{\vs, \vr}^{\omega} \approx F^{\omega} \eta = -\omega^2 \frac{e^{i \omega R}}{\sqrt{R}} \int_{\RR^2}  e^{i \omega ( \mathbf{s} - \vr) \cdot \mathbf{y}} \eta(y) d\vy 
\end{equation}
is the linearized forward operator acting on the perturbation. 

Solving the inverse problem \eqref{eq:fwi} using the linearized operator in \eqref{eq:far_field_pattern} results in an explicit solution given by the normal equation
\begin{equation}\label{eq:FBP}
    \eta^* = \left ( \left(F^{\omega} \right )^* F^{\omega} + \epsilon I \right )^{-1} \left(F^{\omega} \right )^*  \Lambda_{\vs, \vr},
\end{equation}
which is often also referred to as filtered back-projection \cite{Colton_Kress:Integral_Equation_Methods_in_Scattering_Theory}. The constant $\epsilon$ is a small regularization parameter that remedies the ill-conditioning of $\left(F^{\omega} \right )^* F^{\omega}$.
In this case, $\left(F^{\omega} \right )^* F^{\omega}$ is translation invariant thus it is a convolutional operator, and $\left(F^{\omega} \right )^*$ is a Fourier transform thus highly oscillatory. 

\begin{figure}
    \centering
    \includegraphics[width=0.95\textwidth]{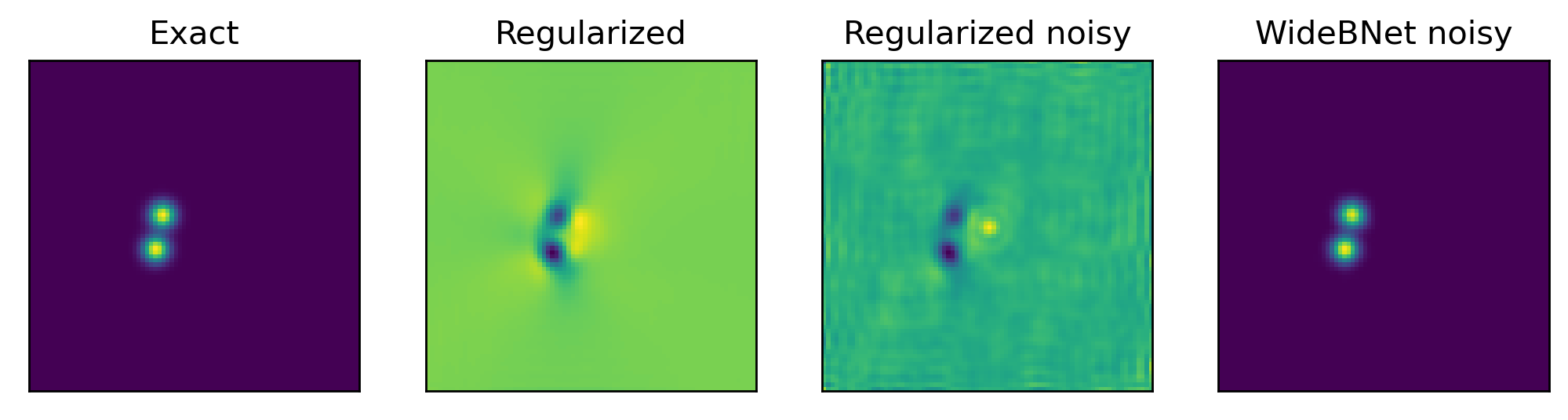}
    \includegraphics[width=0.95\textwidth]{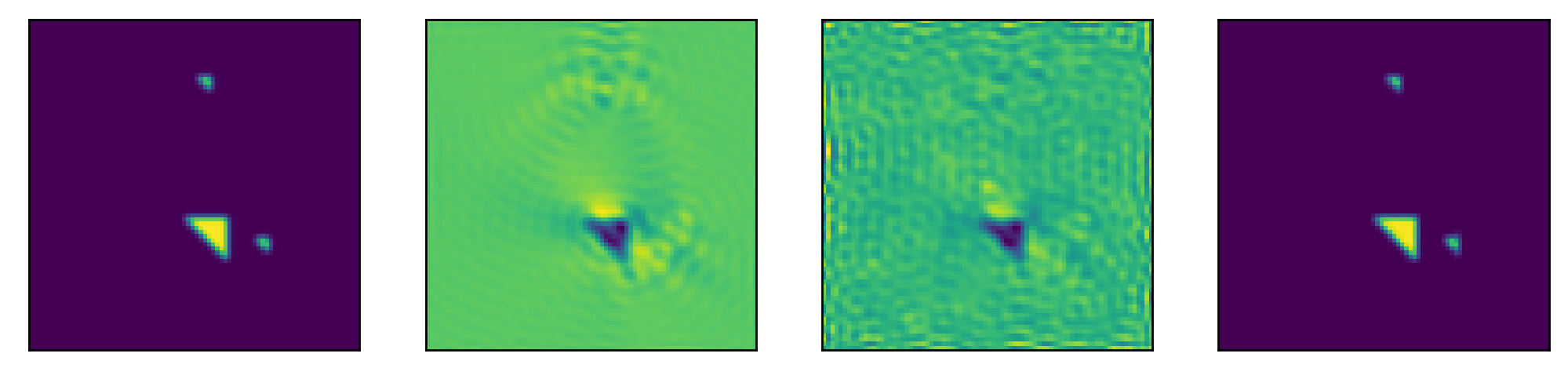}
    \includegraphics[width=0.95\textwidth]{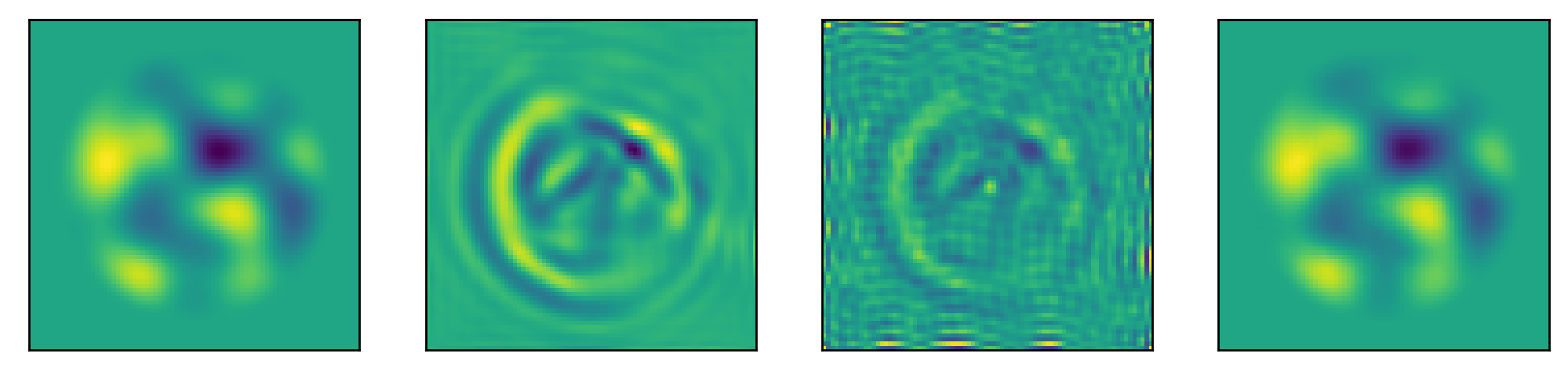}
    \caption{Comparison of reconstruction for inversion of the FIO in \eqref{eq:FBP} using data at $10$ Hz and \wbnn using wideband data at $2.5$, $5$ and $10$ Hz. The two middle images have been rescaled to visually match the other two. The noise is injected to the data following \eqref{eq:noise:model} with $\sigma =1$.}
    \label{fig:FIO_wbnet}
\end{figure}
We point out that in the case that the background slowness squared is non-constant the operator in \eqref{eq:far_field_pattern} becomes a FIO similar to \eqref{eq:Fourier_operator} where $a$ and $\phi$ depend on the background medium, which bring us back to the case covered in \eqref{eq:fwi}.

Finally, we used \eqref{eq:FBP} to solve the inverse scattering problem for the noiseless and noisy case for different perturbations (see Figure \ref{fig:FIO_wbnet}). In this case the regularization parameter $\epsilon$ was chosen, after a labour intensive search, to be $10^{-3}$ and $1$ in the noiseless and noisy cases, respectively. The linear system was solved iteratively using GMRES \cite{Saad_Schultz:GMRES} with a tolerance of $10^{-4}$ and a restart of $10$ steps. We can observe in Figure \ref{fig:FIO_wbnet} that for the perturbations in which the assumptions of weak scattering hold, the filtered back-projection is able to localize the reflectors, albeit with some artifacts that become dominant as noise is incorporated. When the perturbation becomes delocalized or involving several scales the reconstruction deteriorates greatly in both noisy and noiseless cases. In contrast, \wbnn is completely agnostic to these issues and it provides an accurate reconstruction in the three different cases (which the network did not see during the training stage) regardless of the noise.

\section{Implementation Details}
\label{app:architecture_Wide_band}

We provide a succinct and self-contained description of the architecture. For a more detailed exposition please see \cite{MTC_LD_LZN:Wideband_butterfly_net}.

\subsection{Input formatting} \label{sec:input_data}
We assume the scatterers (discretized over an $n_x \times n_z$ grid) and the scattered data (an $n_\textrm{src} \times n_\textrm{rcv}$ matrix for each frequency~$\omega$) are represented using complete quad-trees with $L$~levels
with leaf size~$s$. This translates to a discretization into $n = 2^{\ell}s$ points for each matrix dimension, where $\ell \in [L/2,~L]$, the level, indexes the size of the contiguous $2^{L-\ell}s \times 2^{L-\ell}s$ sub-matrices. 
The choice of $L$ and $s$ is informed by the wavelength scaling in each problem such that the data inside each $s \times s$ cell 
is non-oscillatory, i.e. contains a fixed amount of oscillations. 

\begin{figure}
    \centering
    \includegraphics[trim = 10mm 0mm 350mm 0mm, clip,  height=0.35\textwidth]{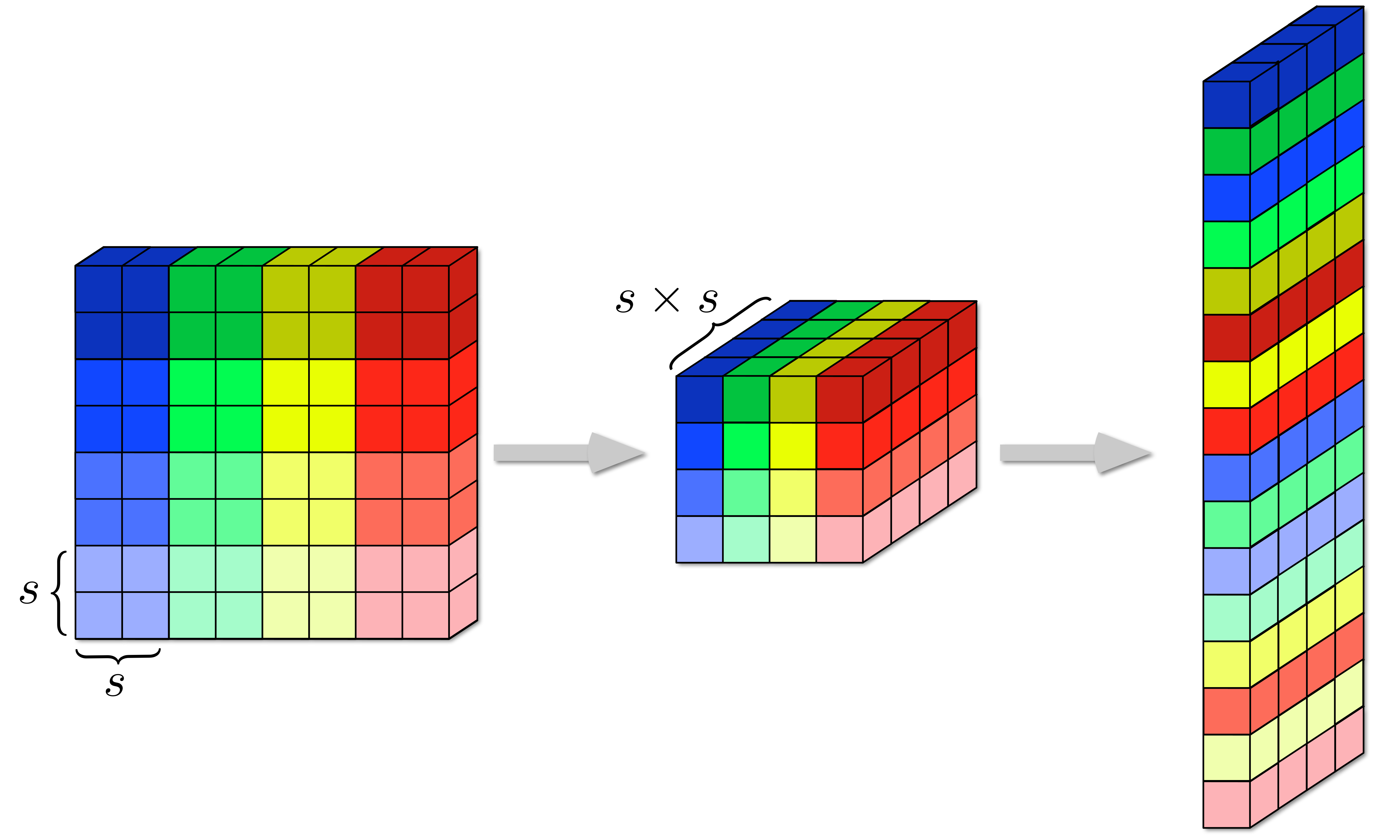}
    \caption{Sketch of the input at each level where the image of dimensions $[2^{\ell}s,~2^{\ell}s]$ is transformed to the tensorized form of size $[2^{\ell},~2^{\ell},~s^2]$, for $\ell = 2$ and $s = 2$.}
    \label{fig:input_ordering}
\end{figure}

Following the Tensorflow convention of \texttt{[height, width, channels]} we reshape these quad-trees into three-tensors of size $[2^{\ell},~2^{\ell},~s^2]$ as shown in Figure~\ref{fig:input_ordering}. The first two dimensions of the tensor contain the geometrical information, and the last contains a local representation. The data describing the local representation inside each voxel corresponds to \emph{channels}. Following this description, we consider the slices along the height and width dimensions, i.e. the geometrical dimensions, as \emph{cells}, e.g. a $1\times1$ cell of data describes slices with dimension $[1,~1,~s^2]$. At the highest spatial resolution \wbnn operates on $1\times1$ cells, and at the lowest spatial resolution it operates on $2^{L/2}\times2^{L/2}$ cells. 

\begin{figure}
    \centering
    \includegraphics[trim = 0mm 0mm 0mm 0mm, clip,  height=0.20\textwidth]{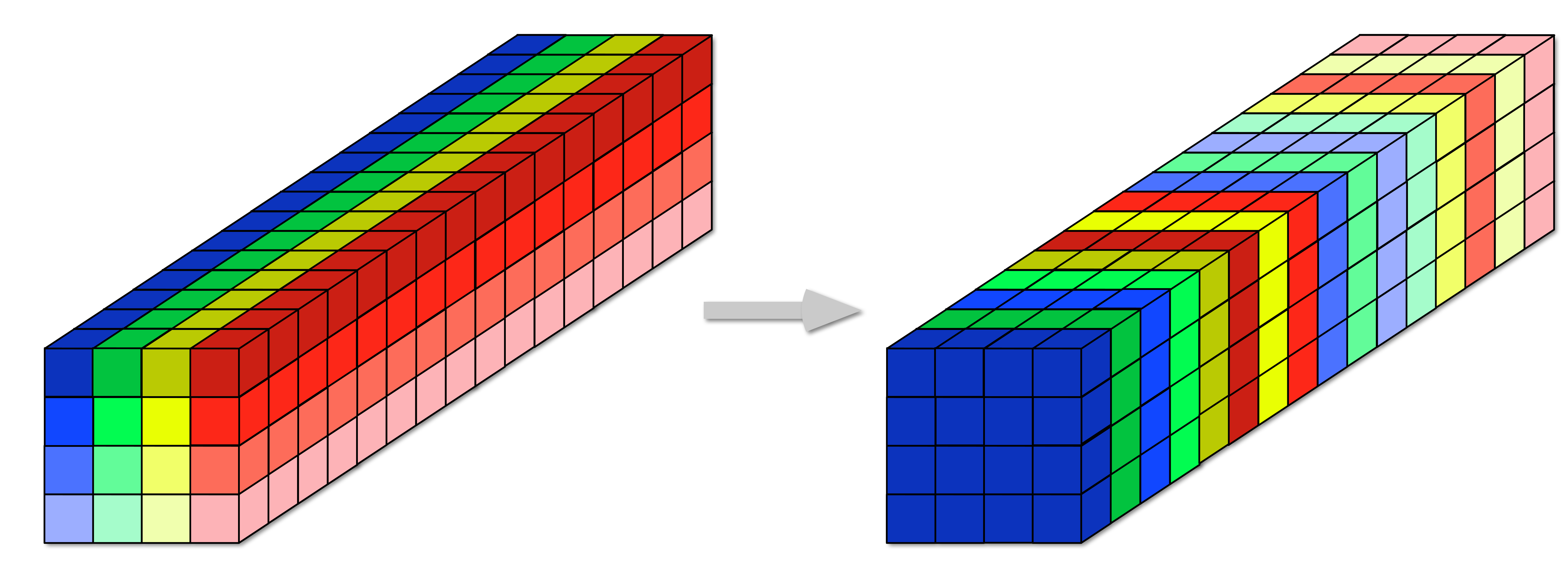}
    \caption{Sketch of the permutation of information of at the $S^{L/2}$ layer. In this case we consider $r=1$, and $L = 4$. The tensor at the middle level is $[2^{L/2}, 2^{L/2}, r 4^{L/2}] = [4,4,16]$. Given that the last dimension is the number of channels, we can observe that the information contained at each cells is distributed along all the other cells. This operation can be easily performed by a combination of reshape and dimension permutation.}
    \label{fig:switch_ordering}
\end{figure}

\begin{figure}[t!]
\centering
\includegraphics[trim = 0mm 50mm 0mm 40mm, clip, height=0.35\textwidth]{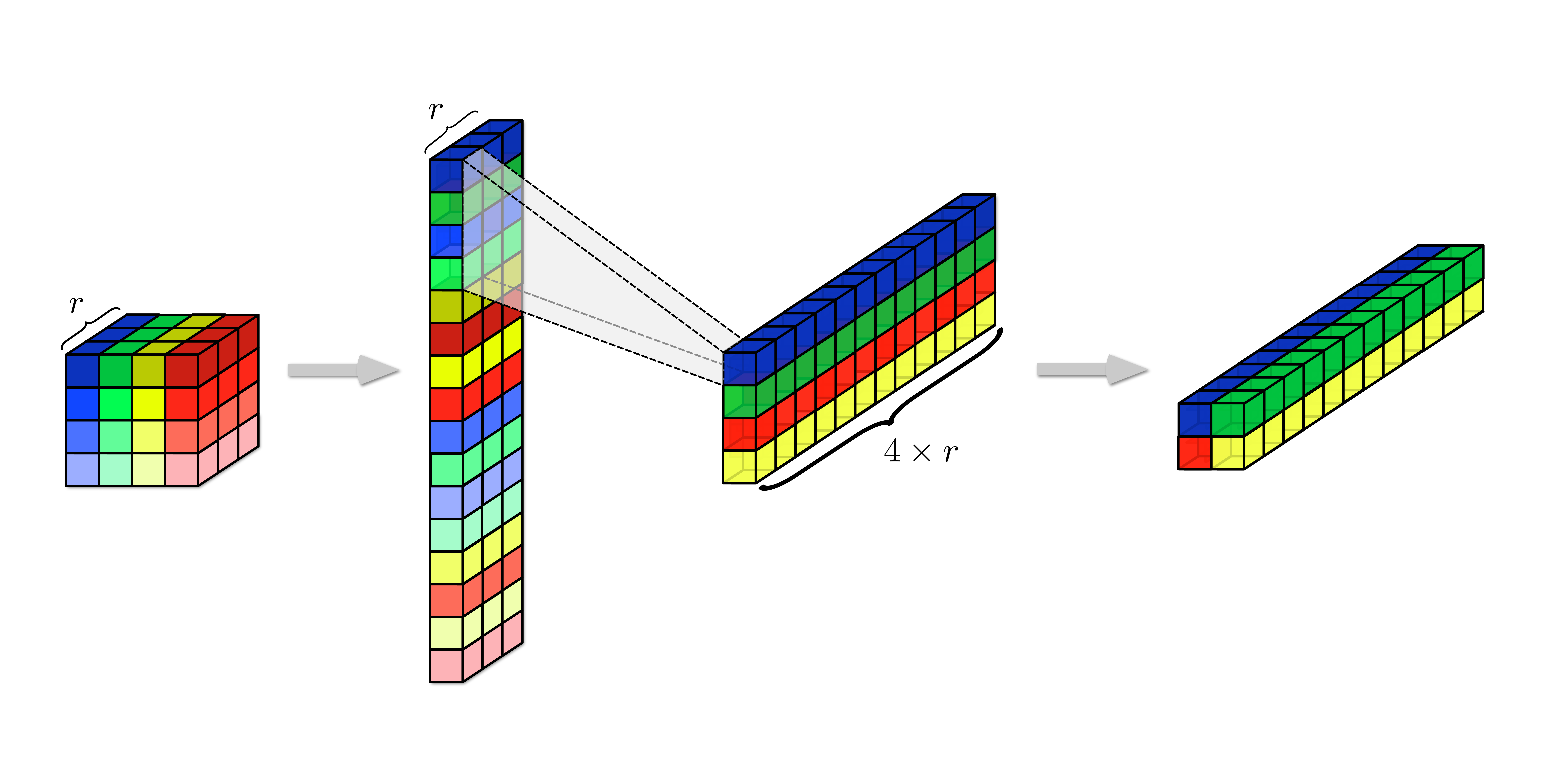}
\caption{
Sketch of the application of the $H^{\ell}$ layer, for convenience we illustrate the reshaping taking place using a Z-ordering. The layer decimates by a factor of four the number of neurons in the spatial dimensions (a factor two in each direction), while increasing four times the number of channels.
\label{fig:H_ell}}
\end{figure}

\begin{figure}[t!]
\centering
    \includegraphics[trim = 0mm 50mm 0mm 40mm, clip, height=0.35\textwidth]{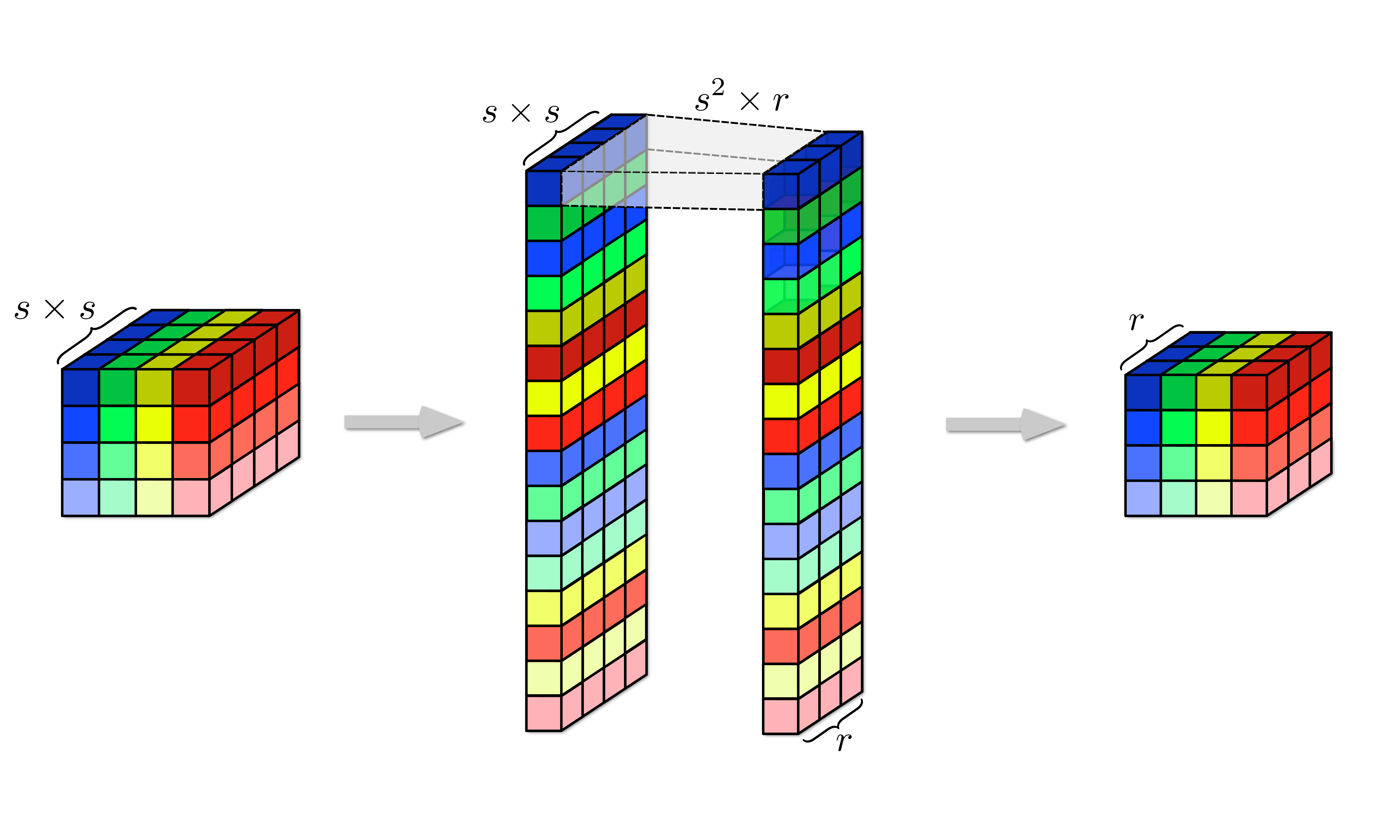}
\caption{Sketch of the compression carried in the $V^{L}$ layer, from the points contained in a leaf of size $s\times s$ to a local representation of rank $r$. We illustrate a reshape operation to a Z-ordering vector to illustrate to local connectivity of the layer. The grey polygon represent the connections between the two layers. \label{fig:V_ell}}
\end{figure}
  
The input data $\Lambda_{\vs,\vr}^{\omega} \in \mathbb{C}^{n_\textrm{\tiny{sr}} \times n_\textrm{\tiny{rs}}}$ is a function of the probe frequency $\omega$ in \eqref{eq:scattering}. 
The input data to be collected from a wide-band of frequencies $ \{ \omega_{\ell}\}_{\ell = L/2}^{L} = \{ \omega/ 2^{L-\ell}\}_{\ell = L/2}^{L}$ 
, i.e, we use a dyadic partition containing $L/2+1$ intervals. With a slight abuse of notation we denote the resulting dataset as $\Lambda_{s,r}^{\ell} : = \Lambda_{s,r}^{\omega/2^{L-\ell}} \in \mathbb{C}^{n_\textrm{src} \times n_\textrm{rcv}}$. 
Following the quad-tree structure we reshape each data tensor $\Lambda_{s,r}^{\ell}$ into a three-tensor of size $[2^{\ell},~2^{\ell}, s^2]$. The input to \wbnn thus consists of the collection $\{\Lambda_{s,r}^{\ell}\}_{L/2 \leq \ell \leq L}$.

\subsection{Architecture Overview} \label{sec:arch_overview}

A schematic diagram of \wbnn is shown in Figure~\ref{fig:WBNet}.
At a high level, \wbnn consists of $L+3$~specialized layers that are non-linear analogues of the butterfly factors in the butterfly factorization. In this case all the layers are linear, and we introduce non-linearities \textsc{resnet} module. These layers ultimately send data into a \textsc{convnet} module, which mimic the effects of the regularized pseudo-inverse  in sharpening the image/estimate as shown in \eqref{eq:FBP}.


The main advantage of \wbnn compared to other butterfly-based architectures, such as \texttt{BNet} \cite{Yingzhou2018} or \texttt{SwitchNet} \cite{Khoo_YingSwitchNet:2019} networks, is in how multi-frequency datasets are assimilated by exploiting the connection between spatial resolution and frequency. In particular, we stress that it is both
\begin{enumerate}[label=(\roman*)]
\item the connectivity/permutations inside the specialized layers $\{H^{\ell}\}$ and $\{G^{\ell}\}$ that process the wide-band data, as well as,
\item the non-linearities induced by the middle $S^{L/2}$ and \textsc{resnet} layer, as shown in Figure \ref{fig:WBNet}, 
\end{enumerate}
that are crucial towards achieving stable training dynamics, as well as image super-resolution, which is argued in \cite{Donoho:super_resolution1992}.

By choosing a dyadic partition of the frequency band \wbnn exploits the inherent multiscale nature of the $\{H^{\ell}\}$ layers. Each $\{H^{\ell}\}$ layer only locally interpolates the data over voxel cells of effectively $2^{L-\ell} \times 2^{L-\ell}$, i.e. the effective length scales at this layer are of order $2^{L-\ell}$. It follows from the dispersion relation in wave-scattering that only data of bandwidth $\omega^{\ell}$ are informative at this length-scale
. 
This strategy of dyadically partitioning the bandwidth to localize spatial information is also employed by the Cooley-Tukey FFT algorithm to achieve quasi-linear time complexity \cite{Cooley_Tukey:1965}; in our setting this strategy affords us significant reductions in the number of trainable weights in the network, thus reducing the number of training samples and the computational complexity.

\subsection{Local embedding and sampling layers} \label{sec:UV_modules}
The $V^{\ell}$ layer can be represented by a block diagonal matrix with block size $r \times s^2$. This layer takes input data (viewed as a complete quad-tree) and compresses the leaf nodes at level~$L$, each with $s \times s$ degrees of freedom, into $1\times1$ cells with $\sqrt{r} \times \sqrt{r}$ degrees of freedom as depicted in Figure~\ref{fig:V_ell}. Similarly, the $U^{L}$ layer can also be represented by block diagonal matrix, albeit with block sizes of $s^2 \times r$. This layer ``samples'' the local representation back to its nominal dimensions. In both instances the compression/decompression is essentially lossless provided the number of levels~$L$ is properly adapted to the probe frequency~$\omega$. Provided that these parameters are chosen correctly, then the data is non-oscillatory (i.e. sub-wavelength) over $s \times s$ length-scales and therefore admits a low-rank representation with rank~$r$.

We also utilize layers $V^{\ell}$ for $L/2 \leq \ell \leq L$ whose inputs are assumed to be sampled from bandwidth $\omega^{\ell}$. Each $V^{\ell}$~layer compresses the input data \emph{at level}~$\ell$ such that nodes with $2^{L-\ell}s \times 2^{L-\ell}s$ degrees of freedom are mapped into $1\times1$ cells with $\sqrt{r} \times \sqrt{r}$ degrees of freedom; this also has the interpretation of spatial downsampling. Note that the dyadic scaling in the definition of $\omega^{\ell}$ is critical in maintaining the balance between spatial resolution and frequency. 

When the input data $\Lambda_{s,r}^{\omega_{\ell}}$ is represented as a two-tensor of dimension $[2^{\ell},~2^{\ell}]$, each $V^{\ell}$ layer can be implemented as a \texttt{LocallyConnected2D} layer in Tensorflow with $r$~channels and both the kernel size and stride as $2^{L-\ell}\times 2^{L-\ell}$. The $U^{L}$ layer can also be implemented as \texttt{LocallyConnected2D} layer with rank~$s^2$ and $1\times1$ kernel size and stride; the input to this layer is assumed to be of dimension $[2^{\ell},~2^{\ell},~c]$ with $c$~input channels.

\subsection{Down- and up-sampling layers} \label{sec:HG_modules}
The $H^{\ell}$ and $G^{\ell}$ layers in Figure \ref{fig:WBNet} continue the theme of multiscale processing. When viewed as matrices, both $H^{\ell}$ and $G^{\ell}$ are block diagonal with block size $4^{L-\ell}r \times 4^{L-\ell}r$. Equivalently, when the input is formatted as a complete quad-tree, this implies both are \textit{local operators} which process the nodes on the tree at length scale $l$ to map each $2^{L-\ell}s \times 2^{L-\ell}s$ cells. Within each block there is further structure to the operators, as Figure \ref{fig:H_ell} demonstrates. For each $H^{\ell}$ each sub-block has the interpretation of \emph{aggregating} information, whereas each $G^{\ell}$ achieves the dual task of \emph{spreading} information. We stress, however, that the action of this is entirely local within each cell. The key observation is that by permuting each node following a set pattern each operator becomes block-diagonal with block size $4^\ell$, for all $L/2 \leq \ell \leq L$.



The $H^{\ell}$ layers differ in that they process two inputs: one the output of the $V^{\ell}$ layer of dimension 
$[2^{L-\ell}, 2^{L-\ell}, r]$, the other the output from the previous layer of dimension $[2^{\ell}, 2^{\ell}, c]$ for some channel size $c$. To process the dimensions of both we first upscale each patch with redundant information to convert the data into $[2^{\ell}, 2^{\ell}, r]$. Then this is concatenated with the other dataset to form a tensor of size $[2^{\ell}, 2^{\ell}, c+r]$.    



\subsection{\textsc{Switch-Resnet} layer} \label{sec:switch_layer}
We retain the permutation pattern of the switch layer as this is responsible for capturing the inherent non-locality of wave scattering (e.g. a point scatterer generates a diffraction pattern that is measured by all receivers in our geometry). We illustrate this pattern in Figure~\ref{fig:switch_ordering}.

The input in this level serves as a condensed representation of the measured data. It is at this level that we non-linearly process the multifrequency dataset; we speculate that this also essential in facilitating the model to produce super-resolved images. We achieve this by using a residual network to refine each channel locally.

\subsection{\wbnn Parameter Count} \label{sec:network_cost}
A tedious but straight-forward computation shows that the total number of parameters  scales as $\mathcal{O} \left ( N(\log N + \log^3 N)\right)$. Note this is essentially linear in the total degrees of freedom in the data ($N$) up to polylogarithmic factors. Furthermore, note if na\"ively $L$ separate single channel \wbnn networks were used to compute \eqref{eq:fwi_weighted} this would correspond to complexity $ \mathcal{O} \left ( N\log N^2 \right)$; the multi-frequency assimilation only exceeds this with mild oversampling by a logarithmic factor.

Lastly, we note the effect of the partitioning of the frequencies. If all the frequencies were ingested at length scale $L$ then the scaling becomes $\mathcal{O}\left (N(\log N^2 + \log N^3) \right)$. While to leading order this presents the same asymptotic scaling, in terms of practical considerations this presents as substantial increase in the number of trainable parameters.

\section{Training details}
\label{app:train_details}

\subsection{Generation of Data}

The scattering data was generated with 5 point-stencil second-order finite differences for training and 9 point-stencil fourth-order finite differences for testing. The radiating boundary conditions  were implemented using PMLs spanning one wavelength with a quadratic profile of intensity~$80$~\cite{Berenger1994}. The domain of interest was discretized using a regular grid of $80 \times 80$ grid points. The scattered waves were sampled at $80$~equi-angular gridpoints in $D$ (see Figure \ref{fig:scattering}) that was implemented as a circle of radius $r=0.5$. The incident waves arrived at angles aligned with the receiver geometry. The training and testing split was $21000$ points and $4000$ points, respectively. For the experiments we choose three frequencies $2.5$, $5$ and $10$Hz, and for simplicity we generated the data in the same mesh, and we use the same sampling geometry for the three frequencies. Thus each datapoint comprised by a tuple of a tensor of dimension~$80 \times 80 \times 3$, and a discretization of the domain of interest $\Omega$ as a matrix of dimension~$80 \times 80$. The data was generated in Matlab using the UMFPACK \cite{Davis:UMFPACK} to solve the linear systems, for each dataset the generation of the data took roughly a week running on a dual-socket Intel Xeon CPU E5-2670 with $384$ GB of RAM. Figure \ref{fig:noisy data} (bottom) depicts one sample of the resulting far-field pattern.  

\subsection{Training}
\begin{figure}
    \centering
    \includegraphics[width=0.80\textwidth]{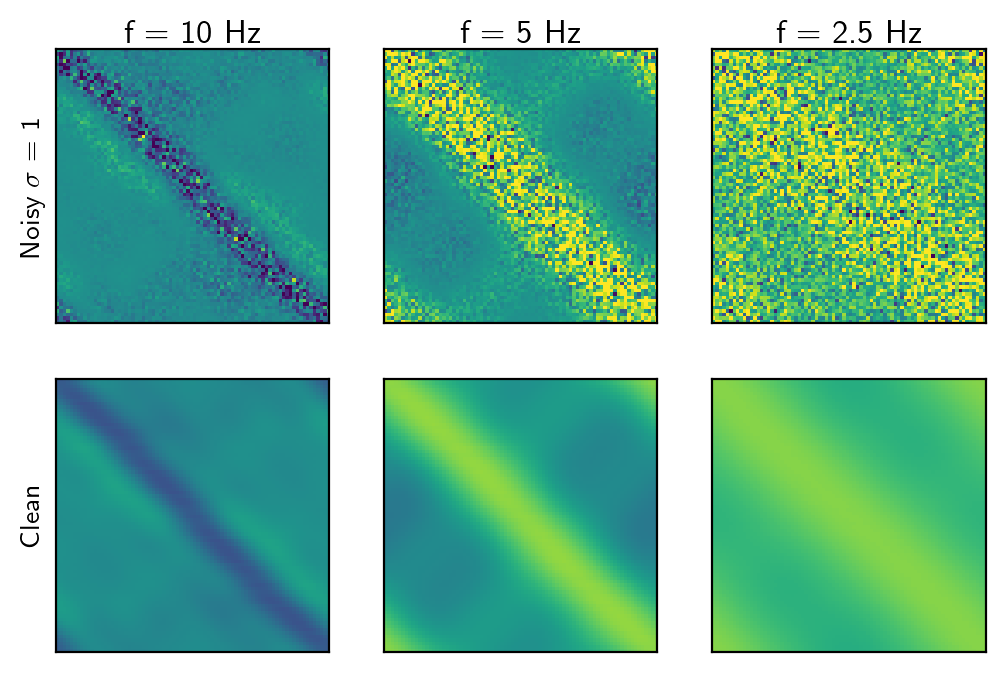}
    \caption{Depiction of effect of $\sigma=1$ multiplicative noise (top) versus clean noiseless data (bottom).}
    \label{fig:noisy data}
\end{figure}

We implemented all the networks in Tensorflow 2.2 and we used the same optimization hyper-parameters for all the networks. The initial learning rate was set to $5e-3$ using an exponential decay schedule with decay rate~$0.95$ every $2000$ plateau steps with staircasing. Optimization was performed by Adam~\cite{kingma2015adam} using the standard parameters $\beta_1 = 0.9$ and $\beta_2 = 0.999$, which we terminate after $501$~epochs. We did not observe any sensitivity to initialization reported in \cite{Butterfly-Net2} and simply used standard glorot uniform sampling for the initial weights. No effort was taken to optimize hyper-parameters on an external validation set. Each experiment was computed with single-precision on four Nvidia GTX-1080Ti graphics cards with 128 GB of shared RAM, each training run was completed in roughly 12 hours. 

The data at each frequency was normalized by subtracting and scaling by the mean and variance of the aggregate pixel intensities in the training dataset. 
We consider a multiplicative coloured noise model similar to \cite{Khoo_YingSwitchNet:2019} with 
\begin{equation}
\Lambda_{\vs,\vr} = \Lambda_{\vs,\vr}^{\text{clean}}(1 + \epsilon(\vs,\vr)),\quad \epsilon(\vs,\vr)\stackrel{iid}{\sim} \mathcal{N}(0, \sigma^2),
\end{equation}
which is applied to each data pixel dynamically each epoch, where we use $\sigma=1$ in order to introduce a $100\%$ noise-to-signal ratio in the training data. Figure~\ref{fig:noisy data} depicts the result of this process.

The loss is the mean square error following 
\begin{equation}
\sum_{x = \text{pixel in image}} \lVert (K_{\text{high}} \ast \eta)(x) - \wbnn[\Lambda_{s,r}^L, \ldots, \Lambda_{s,r}^{L/2}](x)  \rVert_2^2
\end{equation}
where $\eta$ is the total scattering medium and $\{\Lambda_{s,r}^l\}$ the noisy multifrequency data. The scattering medium data was smoothed using a high-pass filter $K_\text{high}$ described by a Gaussian kernel with characteristic width of $0.75$ grid points. This was observed to promote faster training in~\cite{MTC_LD_LZN:Wideband_butterfly_net} for noiseless data, which has the qualitative effect of smoothing the gradient. Since the support of this filter is significantly less than the Nyquist limit this processing does not pollute or eliminate sub-wavelength features in the medium. 

\subsection{Hyper parameters for each architecture} \label{app:hyper_parameters}

\begin{itemize}[leftmargin=*]

    \item  \wbnn: we used $N_{\text{CNN}}=3$ convolutional layers and $N_{\text{RNN}}=3$ residual layers and we used relu activation function throughout the network to add non-linearities. We set the number of levels $L=4$, the leaf size $s = 5$, and the interaction rank $r = 4$. The total number of trainable parameters is $1,914,061$.
    
    \item \nbnn:  we used the same parameters as in \texttt{WideBNet}, but only processing the highest resolution. The total number of trainable parameters is $1,241,293$.

    \item \switchnet: we chose the same parameters as in \cite{khoo2017solving}, i.e.,  rank of the low-rank approximation $t = 3$, the number of partitions in $x$, $P_X = 8^2$, the number of partitions receiver manifold, $P_D = 4^2$, the window size of the convolution layers, $w = 10$, the number of channels in the convolution layers $\alpha = 18$, and the number of convolution layer $L = 3$. For a detailed explanation of these parameters please see  \cite{Khoo_YingSwitchNet:2019}. The total number of trainable parameters is $3,191,416$.

    \item \slbnn: we used the same parameters as in \texttt{WideBNet}. The only difference is the omission of the switch layer $S^{L/2}$, which does not contain any trainable parameters. The total number of trainable parameters is $1,914,061$.
    
    \item \unet: we considered an architecture with five levels, with a width-filter of $3$, we used $4$ channels at the finer level, which was doubled when going to lower levels. To make the implementation more competitive we added a periodic padding in the convolutional layers, in order to maintain the periodicity on the input and output. In addition, we used the full wide-band data, using only narrow-band data provided worse results overall. The total number of trainable parameters is $587,537$.
    
\end{itemize}

\section{Comparison against other architectures}
\label{app:architectures}

\begin{figure}
    \centering
    \includegraphics[width=0.95\textwidth]{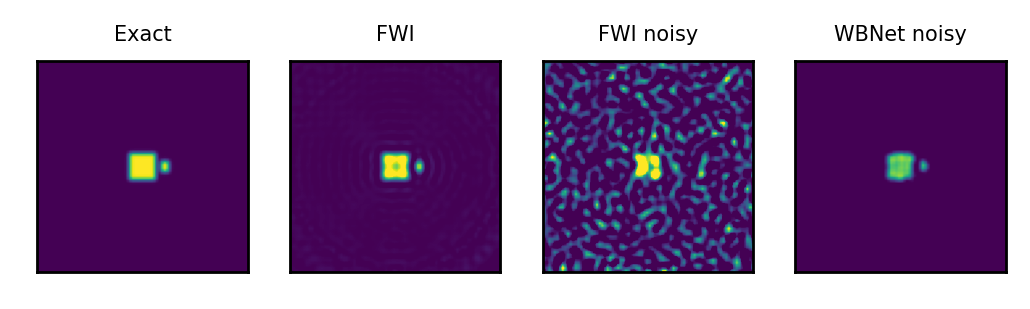}
    \caption{Comparison of reconstruction for the colliding example for $\Delta = \lambda$. From left to right, the reference perturbation to be reconstructed, the reconstruction of the perturbation using full-waveform inversion in Pysit \cite{Pysit} with clean data, FWI with noisy data, and \wbnn with noisy data. The computation time for the noiseless and noisy FWI is around $2$ hours, and $45$ minutes respectively, whereas \wbnn takes $0.1$ seconds for the inference.}
    \label{fig:fwi_bdnet_comparison}
\end{figure}

We compare the resolution with FWI implemented in Pysit \cite{Pysit}(for a comparison against filtered back-projection please see Figure \ref{fig:FIO_wbnet}). We used $20$ iterations to process the data at $2.5$, $5$ and $10$ Hz, in increasing order. The optimization was performed using L-BFGS, plus line-search. We used as an example the collision data-set in which a super wavelenght scatteres is roughly a wavelength away from a sub-wavelength one. In Figure \ref{fig:fwi_bdnet_comparison} we can observe that if the data is clean, FWI is able to recover the perturbation, but when we add data following \ref{fig:noise-sigma}, then the reconstruction becomes much worse. Fortunately, even in this noisy setup \wbnn is able to image the perturbation. 

For completeness we also provide several tables illustrating the superiority of \wbnn with respect to the other architectures considered in this manuscript. Table \ref{tab:rel_l2_training_pixel} provides the mean $\ell^2$ error pixel wise during training. Tables  \ref{tab:rel_l2_training_image} and \ref{tab:rel_l2_testing_image} show the mean relative $\ell^2$ error pixel wise during the training and testing stage for the architectures considered. Tables \ref{tab:rel_l2_training_image} and \ref{tab:rel_l2_testing_image} depict the relative difficulty of learning the different dataset, where we can observe that the blob is, in general, the easiest to learn, whereas, datasets with large back-scattering, such as the Ten-squares and GRF ($\sigma = 0.01)$ are the hardest one to learn.   

\begin{table}[]
\caption{Comparison of mean-squared error per pixel training error between different architectures.}
\label{tab:rel_l2_training_pixel}
\centering
\begin{tabular}{@{}lrrrr@{}}
\toprule
Dataset & \multicolumn{1}{l}{WideBNet} & \multicolumn{1}{l}{UNet} & \multicolumn{1}{l}{SwitchlessBNet} & \multicolumn{1}{l}{NarrowBNet} \\ \midrule
Blob & \textbf{4.548E-04} & 4.270E-03 & 2.364E-03 & 5.185E-04 \\
Gaussian (2h) & \textbf{1.641E-03} & 4.369E-02 & 4.996E-02 & 3.642E-03 \\
Ten Squares (10h) & \textbf{2.642E-02} & 1.051E-01 & 3.479E-01 & 8.022E-02 \\
Right Triangles (3/5/10h) & \textbf{3.347E-03} & 1.209E-01 & 6.175E-02 & 4.186E-03 \\
Shepp Logan & \textbf{6.050E-03} & 3.825E-02 & 2.545E-02 & 9.323E-03 \\
GRF ($\sigma=0.01$) & \textbf{1.082E-02} & 3.924E-02 & 3.135E-01 & 4.579E-02 \\
GRF ($\sigma=0.04$) & \textbf{2.978E-03} & 1.094E-02 & 1.900E-01 & 8.114E-03 \\
Squares (3/5/10h) & 3.412E-03 & 1.037E-01 & 6.731E-02 & \textbf{2.786E-03} \\ \bottomrule
\end{tabular}
\end{table}

\begin{table}[]
\caption{Comparison of relative $\ell^2$ squared error per image in training between different architectures.}
\label{tab:rel_l2_training_image}
\centering
\begin{tabular}{@{}lrrrr@{}}
\toprule
Dataset & \multicolumn{1}{l}{WideBNet} & \multicolumn{1}{l}{UNet} & \multicolumn{1}{l}{SwitchlessBNet} & \multicolumn{1}{l}{NarrowBNet} \\ \midrule
Blob                         & \textbf{4.087E-04} & 4.279E-03 & 2.273E-03 & 5.295E-04 \\
Gaussian (2h)               & \textbf{1.573E-03} & 3.887E-02 & 4.825E-02 & 3.498E-03 \\
Ten Squares (10h)           & \textbf{2.725E-02} & 1.082E-01 & 3.556E-01 & 9.290E-02 \\
Right Triangles (3/5/10h)   & \textbf{3.788E-03} & 1.670E-01 & 7.279E-02 & 6.908E-03 \\
Shepp Logan                 & \textbf{6.091E-03} & 3.851E-02 & 2.565E-02 & 9.471E-03 \\
GRF ($\sigma=0.01$)         & \textbf{1.070E-02} & 3.991E-02 & 3.256E-01 & 6.104E-02 \\
GRF ($\sigma=0.04$)         & \textbf{3.259E-03} & 1.195E-02 & 2.101E-01 & 1.088E-02 \\
Squares (3/5/10h)           & \textbf{3.441E-03} & 1.497E-01 & 8.210E-02 & 5.691E-03 \\ \bottomrule
\end{tabular}
\end{table}

\begin{table}[]
\caption{Comparison of relative $\ell^2$ squared error per image in testing between different architectures.}
\label{tab:rel_l2_testing_image}
\centering
\begin{tabular}{@{}lrrrr@{}}
\toprule
Dataset & \multicolumn{1}{l}{WideBNet} & \multicolumn{1}{l}{UNet} & \multicolumn{1}{l}{SwitchlessBNet} & \multicolumn{1}{l}{NarrowBNet} \\ \midrule
Blob & \textbf{5.449E-04} & 4.493E-03 & 2.598E-03 & 5.955E-04 \\
Gaussian (2h) & \textbf{2.299E-03} & 5.502E-02 & 1.048E-01 & 4.315E-03 \\
Ten Squares (10h) & \textbf{3.074E-02} & 1.079E-01 & 4.060E-01 & 9.098E-02 \\
Right Triangles (3/5/10h) & 7.097E-03 & 1.845E-01 & 1.197E-01 & \textbf{6.662E-03} \\
Shepp Logan & \textbf{7.320E-03} & 3.791E-02 & 3.209E-02 & 1.160E-02 \\
GRF ($\sigma=0.01$) & \textbf{1.482E-02} & 4.811E-02 & 5.665E-01 & 5.964E-02 \\
GRF ($\sigma=0.04$) & \textbf{3.800E-03} & 1.447E-02 & 2.551E-01 & 1.101E-02 \\
Squares (3/5/10h) & 7.085E-03 & 1.889E-01 & 1.528E-01 & \textbf{6.032E-03} \\ \bottomrule
\end{tabular}
\end{table}

\section{Generalization}
\begin{figure}
    \centering
    \includegraphics[width=0.75\textwidth]{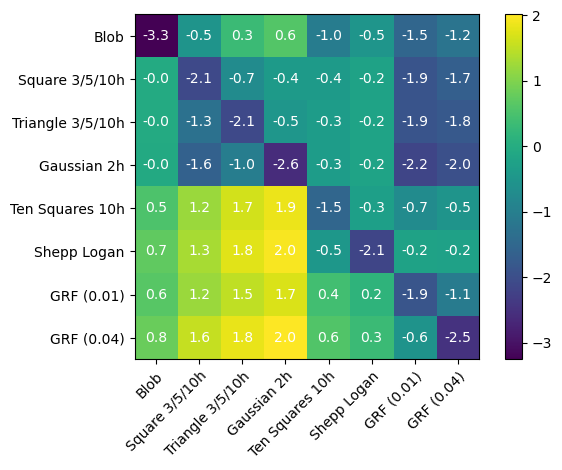}
    \caption{Generalization matrix in which we showcase the average pixelwise squared generalization error for \wbnn models trained with data sets in the ordinate axis, and tested against data sets in the coordinate axis. All numbers are in logarithmic scale (base $10$).}
    \label{fig:generalization_heat_map}
\end{figure}

In Figure~\ref{fig:generalization_heat_map} we provide a quantitative examination of the extrapolative properties of \wbnn models trained and tested with disparate datasets. We observe that training with the square 3/5/10h also results in low testing error on the triangle 3/5/10h and gaussian 2h.

\section{Super Resolution}
\begin{figure}
    \centering
    \includegraphics[width=0.95\textwidth]{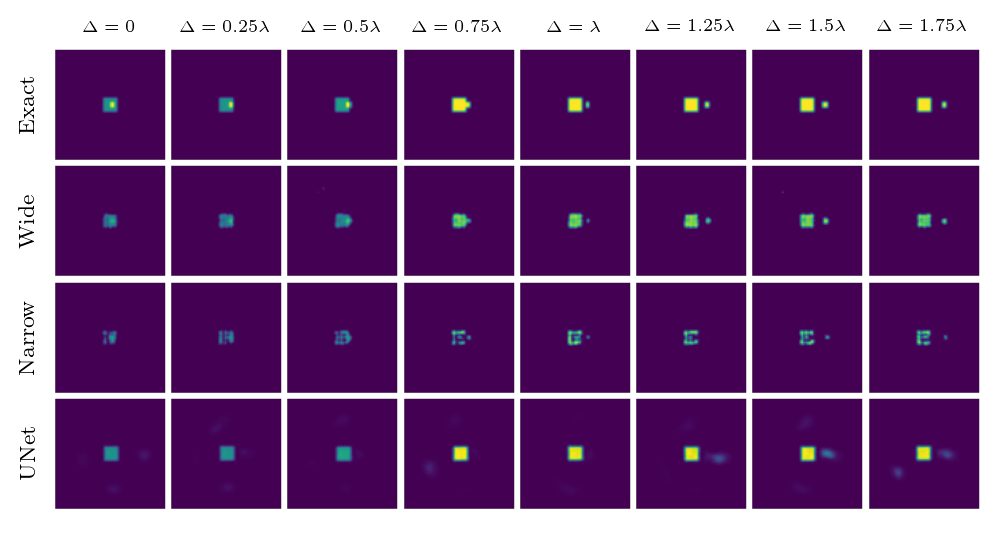}
    \caption{Comparison of imaging capabilities of different networks for square scatterers with side lengths 3 and 10 pixels, which are separated by~$\Delta$. We vary $\Delta$ between sub-wavelength and wavelength separations (one wavelength~$\lambda$ = 8 pixels)}
    \label{fig:collide-square-10-3-h}
\end{figure}

\begin{figure}
    \centering
    \includegraphics[width=0.95\textwidth]{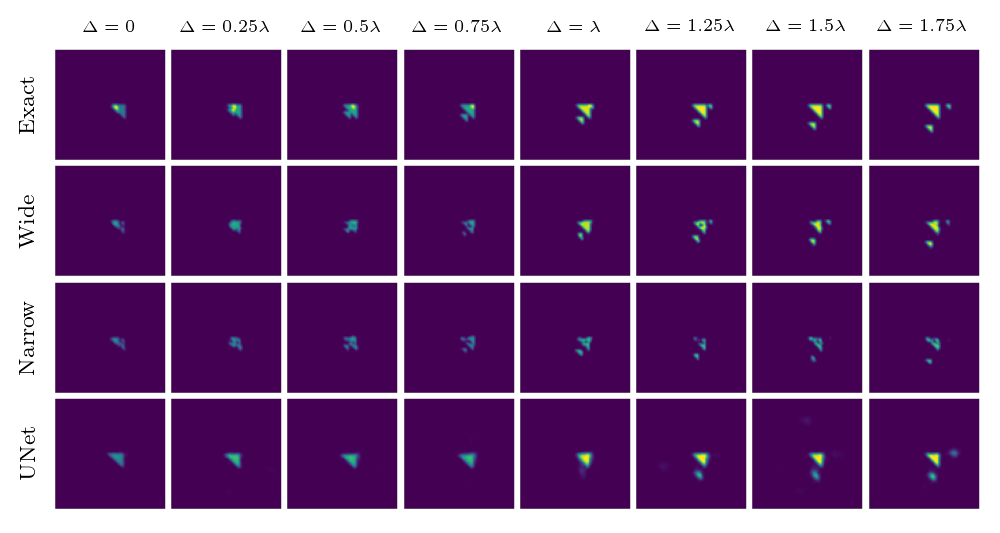}
    \caption{Comparison of imaging capabilities of different networks for triangular scatterers with side lengths 3, 5, and 10 pixels, which are separated by~$\Delta$. We vary $\Delta$ between sub-wavelength and wavelength separations (one wavelength~$\lambda$ = 8 pixels)}
    \label{fig:collide-right-tri-10-3-5-h}
\end{figure}

\begin{figure}
    \centering
    \includegraphics[width=0.95\textwidth]{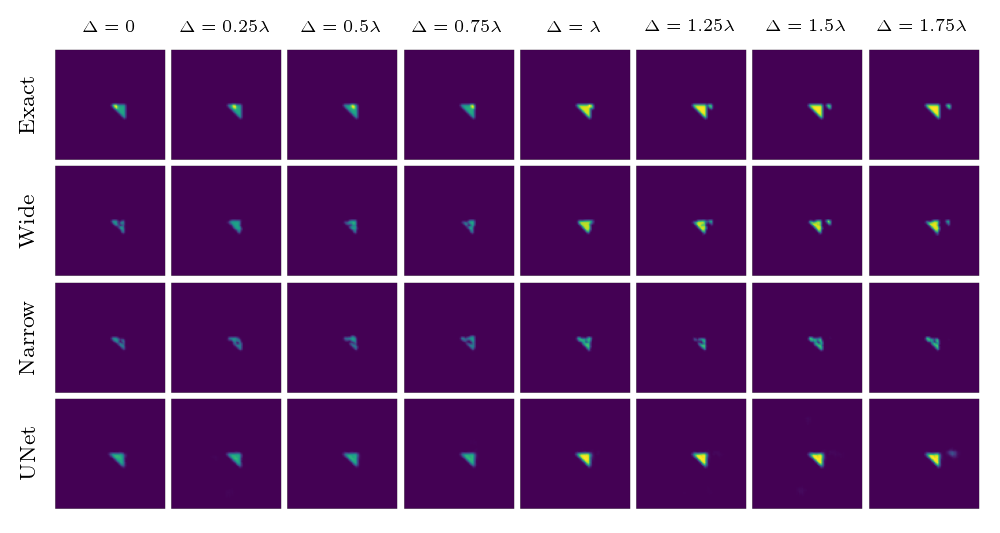}
    \caption{Comparison of imaging capabilities of different networks for triangular scatterers with side lengths 3 and 10 pixels, which are separated by~$\Delta$. We vary $\Delta$ between sub-wavelength and wavelength separations (one wavelength~$\lambda$ = 8 pixels)}
    \label{fig:collide-right-tri-10-3-h}
\end{figure}

For completeness we also provide results for super-resolution beneath the diffraction limit for two square scatterers (Figure~\ref{fig:collide-square-10-3-h}), three triangles (Figure~\ref{fig:collide-right-tri-10-3-5-h}), and two triangles (Figure~\ref{fig:collide-right-tri-10-3-h}). Similar observations discussed in Section~\ref{ssec:superres} hold for these examples as well -- both the \nbnn and \unet architectures struggle in the subwavelength ($\Delta < \lambda$) regimes, whereas \wbnn is still able to consistent image the scatterers.

\end{document}